\documentclass[leqno,12pt]{article}
\begin{document}
\def\gint#1{{\bf [}#1{\bf ]}}
\def\A{{\rm{\bf A}}}
\def\modd#1{\vert#1\vert}
\def\Hs{\heartsuit}
\def\div{{\rm div}}
\def\Arr{{\cal R}}
\def\nrm#1{\Vert#1\Vert}
\def\dq{\lq\lq}
\def\fstar{f^*}
\def\quote#1{\lq\lq#1''}
\def\Rd{{{\bf R}^d}}
\def\Upperhalfspace{{{\bf R}^{d+1}_+}}
\def\Dee{{\cal D}}
\def\Bee{{\cal B}}
\def\Cinf{{\cal C}^\infty_0(\Rd)}
\def\Lwunloc{L^1_{loc}(\Rd)}
\def\Lwunlocmu{L^1_{loc}(\mu)}
\def\Lwunlocv{L^1_{loc}(v)}
\def\fn#1#2{\footnote{#1}{#2}}
\def\fnno#1#2{\footnote{$ ^{#1}$}{#2}}
\def\sqrterm{{\nrm{a_Q(f)}^2_2\over\modd Q}}
\def\Ainf{A_\infty}
\def\Gee{{\cal G}}
\def\Eff{{\cal F}}
\def\M{{\cal M}}
\def\W{{\cal W}}
\def\I{{\cal I}}
\def\Eee{{\cal E}}
\def\Sss{{\cal S}}
\def\N{{\bf N}}
\def\cs{{$\clubsuit$}}
\def\Uball{\{x:\ \modd x\leq1\}}
\def\Ca{Calder\'on}
\def\Pf{{\bf Proof.}\ \ }
\def\AN{{\cal A}_N}
\def\Ckay{{\cal C}^k(\Rd)}
\def\Sf#1#2#3{\left(\int\limits_{\modd{x-t}<#3 y}\modd{#1 * #2_y(t)}^2\,
{dt\, dy\over y^{d+1}}\right)^{1/2}}
\def\Sq#1#2#3{S_{#2,#3}(#1)}
\def\R{{\bf R}}
\def\bigmodd#1{\left\vert#1\right\vert}
\def\bignrm#1{\left\Vert#1\right\Vert}
\def\stars{\smallskip*\hfil*\hfil*\smallskip}
\def\Bispace{{\R^{d_1}\times\R^{d_2}}}
\def\leaderfill{\leaders\hbox to 1em{\hss.\hss}\hfil}
\def\topline#1{\headline={\ifnum\pageno>1 wilson\hfil #1\else\hfil\fi}}
\def\bigbrace#1{{\left\lbrace#1\right\rbrace}}
\def\Calpha{{{\cal C}_\alpha}}
\def\CalphaM{{{\cal C}_{\alpha,M}}}
\def\CalphaO{{{\cal C}_{\alpha,0}}}
\def\CalphaMO{{{\cal C}_{\alpha,0,M}}}
\def\qed{{$\clubsuit$}}
\def\Compnum{{\bf C}}
\def\rmk{\emph{Remark}}
\def\C{\Compnum}
\newtheorem{theorem}{Theorem}
\newtheorem{lemma}{Lemma}
\newtheorem{definition}{Definition}
\newtheorem{corollary}{Corollary}

\title{Stability of Haar decompositions}\author{Michael Wilson\date{}}

\maketitle

\begin{abstract}
We prove a general result implying the $L^2$ stability of Haar decompositions of $L^2(\Rd)$ functions when the Haar functions are distorted by arbitrary, independent, affine changes of variable that are close to the identity. We apply our method to get fully $d$-dimensional generalizations of results of Aimar, Bernardis, Gorosito, Govil, and Zalik, on constructing frames of smooth functions which are, in many natural senses, arbitrarily close to the Haar functions. We also obtain a best-possible estimate on the $L^2$ sensitivity of dyadic averages of functions to small distortions caused by local affine changes of variable.
\end{abstract}

\section{Introduction}

A set $\{\psi_\gamma\}_{\gamma\in\Gamma}\subset L^2(\Rd)$ is called a \emph{Bessel family} if there is a $0\leq B<\infty$ so that
\begin{equation}\label{besselfambound}\sum_{\gamma\in\Gamma}\modd{\langle f,\psi_\gamma\rangle}^2\leq B\nrm f_2^2\end{equation}
for all $f\in L^2(\Rd)$, where $\langle\cdot,\cdot\rangle$ means the usual $L^2$ inner product: $\langle f,g\rangle:=\int f\,\overline g\,dx$. The minimal such $B$ is called the family's optimal bound. By duality, (\ref{besselfambound}) holds for all $f\in L^2(\Rd)$ if and only if, for all finite $\Eff\subset\Gamma$ and all $\{\lambda_\gamma\}_{\gamma\in\Eff}\subset\C$,
$$\bignrm{\sum_\Eff\lambda_\gamma\psi_\gamma}_2^2\leq B\sum_\Eff\modd{\lambda_k}^2.$$
This latter condition is often called \emph{almost-orthogonality}. 

The separability of $L^2(\Rd)$ implies that Bessel families are countable (can have only countably many non-zero elements), so that any Bessel family can be written as a \emph{Bessel sequence} indexed over a subset of \N\ (usually all of \N). Families are preferable to sequences when it's convenient to index the functions over the dyadic cubes in $\Rd$ (the most important case for this paper).

For any non-empty set $\Gamma$, we can turn the collection of Bessel families indexed over $\Gamma$ into a vector space by defining the basic operations in the obvious ways:
\begin{eqnarray*}\{\psi_\gamma\}_{\gamma\in\Gamma}+\{\phi_\gamma\}_{\gamma\in\Gamma}&:=&\{\psi_\gamma+\phi_\gamma\}_{\gamma\in\Gamma}\\
\alpha\{\psi_\gamma\}_{\gamma\in\Gamma}&:=&\{\alpha\psi_\gamma\}_{\gamma\in\Gamma}.\end{eqnarray*}
We call the resulting space $AO(\Gamma)$. The {square root} of the optimal Bessel bound for $\{\psi_\gamma\}_{\gamma\in\Gamma}\in AO(\Gamma)$ defines a norm on $AO(\Gamma)$, which we will refer to as 
$$\bignrm{\{\psi_\gamma\}_{\gamma\in\Gamma}}_{AO(\Gamma)}.$$
It is an easy exercise to show that $(AO(\Gamma),\nrm{\cdot}_{AO(\Gamma)})$ is a Banach space. 

Bessel bounds help to quantify how close some function families come to having useful properties, such as being an orthonormal basis for $L^2$. If $\{\psi_\gamma\}_{\gamma\in\Gamma}\subset L^2(\Rd)$ is a complete orthonormal set, then, for every $f\in L^2(\Rd)$,
$$\sum_{\gamma\in\Gamma}\modd{\langle f,\psi_\gamma\rangle}^2=\nrm f_2^2$$
and 
$$f=\sum_{\gamma\in\Gamma}\langle f,\psi_\gamma\rangle \psi_\gamma,$$
with the series converging unconditionally in $L^2$. If $\{\widetilde\psi_\gamma\}_{\gamma\in\Gamma}\subset L^2(\Rd)$ is another family, such that 
$$\bignrm{\{\psi_\gamma-\widetilde\psi_\gamma\}_{\gamma\in\Gamma}}_{AO(\Gamma)}\leq\delta<1,$$ 
then
$$(1-\delta)^2\nrm f_2^2\leq \sum_{\gamma\in\Gamma}\modd{\langle f,\widetilde\psi_\gamma\rangle}^2\leq(1+\delta)^2\nrm f_2^2$$
for any $f\in L^2$, which says that $\{\widetilde\psi_\gamma\}_{\gamma\in\Gamma}$ is a \emph{frame} \cite{ChrFRB}, with lower and upper frame bounds $(1\mp\delta)^2$; and the series
$$\sum_{\gamma\in\Gamma}\langle f,\widetilde\psi_\gamma\rangle\widetilde\psi_\gamma$$
converges unconditionally to some $\widetilde f\in L^2(\Rd)$  (\cite{ChrFRB}; see Corollary 3.2.5), with $\nrm{f-\widetilde f}_2\leq\delta(2+\delta)\nrm{f}_2$. If $\delta$ is small then, for many signal processing purposes, the family $\{\widetilde\psi_\gamma\}_{\gamma\in\Gamma}$ is almost as good as $\{\psi_\gamma\}_{\gamma\in\Gamma}$. 

The Haar functions might be the simplest example of a complete orthonormal set in $L^2(\Rd)$. They are defined this way. An interval $I\subset\R$ is called dyadic if it has the form $[j2^n,(j+1)2^n)$ for integers $j$ and $n$. We denote the set of dyadic intervals by $\Dee_1$. If $I\in\Dee_1$ we define
\begin{equation}\label{haarfunctiondef}h^{(I)}(x):=\cases{1&if $x\in I_l$;\cr -1&if $x\in I_r$;\cr0&if $x\notin I$;\cr}\end{equation}
where $I_l$ and $I_r$ are $I$'s left and right halves. We use $\ell(I)$ to mean $I$'s length and $x_I$ to mean its geometric center. 

The Haar functions on \R\ are the family $\{h^{(I)}/\modd I^{1/2}\}_{I\in\Dee_1}$, where, here and in the sequel, $\modd E$ means $E$'s Lebesgue measure when $E$ is a measurable set. (The dimension of the measure will change; context will make it clear.)
They are orthonormal; a little more work shows them to be complete. If we set $h_{(I)}:=h^{(I)}/\modd I^{1/2}$ then it's easy to see that the family $\{h_{(I)}\}_{I\in\Dee_1}$ is derivable from $h_{[0,1)}$ by a simple formula: if $I=[j2^n,(j+1)2^n)$ then 
\begin{equation}\label{simpleformula1}h_{(I)}(x)=2^{-n/2}h_{[0,1)}(2^{-n}x-j).\end{equation}
Formulas like (\ref{simpleformula1}) are famililar in frame theory \cite{ChrFRB}.

The Haar functions are easily generalized to $\Rd$. A \emph{cube} $Q\subset\Rd$ is a cartesian product of $d$ intervals of the form $[a,b)$, all having the same length (which we call $Q$'s sidelength and denote by $\ell(Q)$): $Q:=\prod_1^d I_i(Q)$, where $I_i(Q)$ is $Q$'s $x_i$ cartesian factor. We write the geometric center of $Q$ as $x_Q$. We call $Q$ dyadic if all of its factors are dyadic. We denote the set of dyadic cubes in $\Rd$ by $\Dee_d$.

For each $Q=\prod_1^d I_i(Q)\in\Dee_d$ we can define $2^d-1$ functions $h^{(Q),k}$ ($1\leq k<2^d$) via
$$h^{(Q),k}(x)=h^{(Q),k}(x_1,x_2,\ldots,x_d):=\prod_1^d g_i(x_i),$$
where each $g_i$ can be $h^{(I_i(Q))}(x_i)$ or $\chi_{I_i(Q)}(x_i)$, but we don't allow the product $\prod_1^d \chi_{I_i(Q)}(x_i)=\chi_Q(x)$. Then 
\begin{equation}\label{multidimhaar}\bigbrace{{h^{(Q),k}\over\modd{Q}^{1/2}}}_{Q\in\Dee_d\atop 1\leq k<2^d}\end{equation}
is a complete orthonormal system for $L^2(\Rd)$. We can easily make the definitions of $h^{(Q),k}$, for the same $k$, \dq consistent'' across the different $Q$s. Set $Q_0:=[0,1)^d$ and fix some enumeration of the functions $\{h^{(Q_0),k}\}_{1\leq k<2^d}$. If we put $h_{(Q),k}:=h^{(Q),k}/\modd{Q}^{1/2}$ then a formula analogous to (\ref{simpleformula1}) gives all the $h_{(Q),k}$s (and thus all the $h^{(Q),k}$s) in terms of the functions $h_{(Q_0),k}$ ($1\leq k<2^d$). For $Q=\prod_1^d [j_i2^n,(j_i+1)2^n)$, we define
\begin{eqnarray*}&h_{(Q),k}(x_1,x_2,x_3,\ldots,x_d)\\
&:=2^{-nd/2}h_{(Q_0),k}(2^{-n}x_1-j_1,2^{-n}x_2-j_2,2^{-n}x_3-j_3,\ldots,2^{-n}x_d-j_d),\end{eqnarray*}
which obviously gives back the family (\ref{multidimhaar}). 

The Haar functions have a nice, compact form, but their discontinuities give them bad localization in frequency space (i.e., their Fourier transforms decay rather slowly).  
In \cite{GZ} N. K. Govil and R. A. Zalik showed that, for every $\epsilon>0$ and every $m\in\N$, one can create a family of functions $\{\phi_{(I)}\}_{I\in\Dee_1}$ such that 
\begin{equation}\label{govilzalik}\bignrm{\bigbrace{h_{(I)}-\phi_{(I)}}_{I\in\Dee_1}}_{AO(\Dee_1)}<\epsilon\end{equation}
and

a) every $\phi_{(I)}$ has support contained in $(1+\epsilon)I$, where this means the $(1+\epsilon)$-fold expansion centered at $x_I$;

b) every $\phi_{(I)}$ belongs to ${\cal C}^m(\R)$;

c) every $\phi_{(I)}$ equals $h_{(I)}$ at points further than $\epsilon\ell(I)$ away from the discontinuities of $h_{(I)}$; 

d) If $I=[j2^n,(j+1)2^n)$ then $\phi_{(I)}(x)=2^{-n/2}\phi_{[0,1)}(2^{-n}x-j)$.

The proof in \cite{GZ} involved some delicate constructions and estimates with $B$-splines. Later, in \cite{ABG}, H. A. Aimar, A. L. Bernardis, and O. P. Gorosito showed how to obtain families $\{\phi_{(I)}\}_{I\in\Dee_1}$ satisfying (\ref{govilzalik}) and a)-d) by convolving the Haar functions with suitable dilates of a smooth, even function $\psi$. Their work was considered later by M. Chung in \cite{Chung}. In \cite{Z} Zalik extended the result of \cite{GZ} to $d$ dimensions, essentially by taking tensor products of the functions constructed in \cite{GZ} (see \cite{Z}, Lemma 2.5).

The first theorem we prove in this paper is a fully $d$-dimensional extension of the results from \cite{GZ}\cite{ABG}\cite{Z}. We remind the reader that, for $\psi:\Rd\to\C$ and $\delta>0$, $\psi_\delta(x):=\delta^{-d}\psi(x/\delta)$. Here is the theorem.

\begin{theorem}\label{genofabg}Set $\Gamma:=\{(Q,k):\ Q\in\Dee_d,\ 1\leq k<2^d\}$. There are positive constants $C(d)$ and $c(d)$, depending only on $d$, so that the following is true. Let $\psi\in L^1(\Rd)$ be non-negative, with support contained in $[-1,1]^d$ and $\int\psi\,dx=1$. For arbitrary but fixed $0<\eta<1/2$, and for every $Q\in\Dee_d$ and $1\leq k<2^d$, let us define
$$\phi_{(Q),k}(x):=(h_{(Q),k}*\psi_{\eta\ell(Q)})(x).$$
Then
\begin{equation}\label{abgframe}\bignrm{\bigbrace{h_{(Q),k}-\phi_{(Q),k}}_{(Q,k)\in\Gamma}}_{AO(\Gamma)}<C(d)\eta^{1/2}\end{equation}
and the family $\{\phi_{(Q),k}\}_{(Q,k)\in\Gamma}$ satisfies: a) supp $\phi_{(Q),k}\subset (1+2\eta)Q$; b) $\phi_{(Q),k}=h_{(Q),k}$ at points $>c(d)\eta\ell(Q)$ away from $h_{(Q),k}$'s discontinuities; c) if $Q=\prod_1^d [j_i2^n,(j_i+1)2^n)$ then
\begin{eqnarray}\label{genfamily}&\phi_{(Q),k}(x_1,x_2,x_3,\ldots,x_d)\nonumber\\
&:=2^{-nd/2}\phi_{(Q_0),k}(2^{-n}x_1-j_1,2^{-n}x_2-j_2,2^{-n}x_3-j_3,\ldots,2^{-n}x_d-j_d).\end{eqnarray}
\end{theorem}

The convolving function $\psi$ can be as rough or as smooth---and have as much or as little symmetry---as we like. In a sense we'll make precise (which the proof will show), the $\psi$s don't even have to be functions and, if we're willing to give up (\ref{genfamily}), we can choose a different $\psi$ for each $(Q,k)\in\Gamma$.

Our result comes from taking a point of departure different from that in \cite{GZ}\cite{ABG}\cite{Z}. Instead of looking at \dq perturbations'' in the sense of carefully constructed approximations to the Haar functions, we consider \dq perturbations'' in the sense of small errors in dilation and translation. We show that, if the perturbations are small, the resulting families are close to the Haar functions in the $\nrm{\cdot}_{AO(\Gamma)}$ sense. Theorem \ref{genofabg} then follows by averaging over translations.

Our work depends on a sufficient condition for almost-orthogonality (or being a Bessel family) which is based on \emph{bounded variation}. This condition's first (one-dimensional) appearance seems to be in Lemma 2 of \cite{GZ}, which we restate here as Theorem \ref{onedimbvresult}.

\begin{theorem}\label{onedimbvresult}There is a finite $B$ so that the following holds. Let $\{f^{(I)}\}_{I\in\Dee_1}$ be a family of functions $f^{(I)}:\R\to\C$ such that, for every $I\in\Dee_1$: a) $f^{(I)}$ is of bounded variation on \R, with total variation $=:V_f(\R)\leq1$; b) $\{x\in\R:\ f^{(I)}(x)\not=0\}\subset I\}$; c) $\int f^{(I)}\,dx=0$. Then $\{f^{(I)}/\modd I^{1/2}\}_{I\in\Dee_1}$ is a Bessel family, with a Bessel bound $B$. \end{theorem}

What is striking about Theorem \ref{onedimbvresult} is that it does not ask the $f^{(I)}$s to have any obvious smoothness (such as H\"older continuity) or to have all of their discontinuities on carefully placed sets (as with the Haar functions). 

In \cite{WilsonBvconvao} we proved a $d$-dimensional version of Theorem \ref{onedimbvresult}. Its statement requires a definition we will use later.

\begin{definition}\label{nbv}If $f:\Rd\to{\bf C}$ we say that $f\in NBV(\Rd)$\footnote{\dq Normalized bounded variation on $\Rd$''.} if $f$ is measurable and, for every $y\in\Rd$ and every standard coordinate vector ${\bf e}_k$ ($1\leq k\leq d$),  the function $\phi_{{\bf e}_k,y}:\R\to{\bf C}$ defined by
$$\phi_{{\bf e}_k,y}(t):= f(y+t{\bf e}_k)$$
is of bounded variation on $\R$, with $V_{\phi_{{\bf e}_k,y}}(\R)\leq 1$. If $Q\subset\Rd$ is a (possibly non-dyadic) cube we say that $f\in NBV(Q)$ if $f\in NBV(\Rd)$ and $\{x:\ f(x)\not=0\}\subset Q$. If we also have $\int f\,dx=0$ we say $f\in NBV_0(Q)$. 
If $\Gee$ is a set of (possibly non-dyadic) cubes in $\Rd$ and $\{f^{(Q)}\}_{Q\in\Gee}$ is a family of functions such that every $f^{(Q)}\in NBV(Q)$ (respectively $NBV_0(Q)$) we call $\{f^{(Q)}\}_{Q\in\Gee}$ an $NBV$ ($NBV_0$) family. 
\end{definition}

\emph{Remark.} \dq$f\in NBV$'' means $f$ is of bounded variation along every line parallel to one of the coordinate axes, with total variation always $\leq 1$. \medskip

The $d$-dimensional extension of Theorem \ref{onedimbvresult}, proved in \cite{WilsonBvconvao} (Theorem 1), is

\begin{theorem}\label{bvconvao}Let $\{f^{(Q)}\}_{Q\in\Dee_d}$ be an $NBV_0$ family. Then 
\begin{equation}\label{family}\bigbrace{{f^{(Q)}\over\modd{Q}^{1/2}}}_{Q\in\Dee_d}\end{equation}
is a Bessel family, with a bound $\leq ((1+{1\over\sqrt 2})d)^2$.\end{theorem}

\rmk. In \cite{WilsonBvconvao} the result is stated in terms of the almost-orthogonality constant.\medskip

We get Theorem \ref{genofabg} from the following theorem, proved in \cite{WilsonBvconvao} (see Theorem 3 and its proof, especially inequality (21); see also inequalities (30)--(32) and the argument following them).

\begin{theorem}\label{wheeden}There is a constant $C(d)$ so that the following holds. Suppose that $\{f^{(Q)}\}_{Q\in\Dee_d}$ is an $NBV_0$ family and, for every $Q\in\Dee_d$, we have a real $d\times d$ diagonal matrix $D^{(Q)}$ and a vector $y^{(Q)}\in\Rd$ such that
$$\nrm{D^{(Q)}-I_d}_\infty+\nrm{y^{(Q)}}_\infty\leq \eta ,$$
where $\nrm{\cdot}_\infty$ means the naive \lq maximum value' vector norm, $I_d$ is the $d\times d$ identity, and $\eta<1/2$ is a fixed small positive number. Define
$$\widetilde{f^{(Q)}}(x):=f^{(Q)}(D^{(Q)}(x-x_Q+\ell(Q)y^{(Q)})+x_Q).$$
Then
$$\bigbrace{{f^{(Q)}-\widetilde{f^{(Q)}}\over\modd{Q}^{1/2}}}_{Q\in\Dee_d}$$
is a Bessel family, with a Bessel bound $\leq C(d)\eta$.\end{theorem}

\rmk. The result in \cite{WilsonBvconvao} is also stated in terms of the almost-orthogonality constant. Easy examples (see \cite{WilsonBvconvao}) show that the exponent on $\eta$ is sharp: it can't be bumped up to $1+\delta$ for any $\delta>0$. See also Corollary \ref{covidcordyadav} at the end of the paper. \medskip

\rmk. In \cite{WilsonBvconvao} the perturbing changes of variable are given as multiparameter dilations $(x_1,x_2,x_3,\ldots,x_d)\to(\delta_1x_1,\delta_2x_2,\delta_3x_3,\ldots,\delta_dx_d)$ (with $\max_i\modd{1-\delta_i}\leq\eta$), followed by small translations, which is equivalent to the form given in Theorem \ref{wheeden}. In \cite{WilsonBvconvao} the $\ell(Q)y^{(Q)}$s are, in effect, \dq outside'' the parentheses surrounding $x-x_Q$. Experience has shown that putting them \dq inside'' works better, at the cost of a small change in the final constant $C(d)$.\medskip

Every $h^{(Q),k}$ equals 4 times a function in $NBV_0(Q)$. With Theorem \ref{wheeden} in hand, the averaging argument we alluded to above gives Theorem \ref{genofabg} in a few easy steps. 

Trying to extend the results of \cite{GZ}\cite{ABG}\cite{Z} was {not} what led us to explore $NBV$ families and their perturbations in \cite{WilsonBvconvao}. We were interested in how \dq discretization'' (replacing basis or frame functions by averages over small cubes), combined with errors in dilation and translation, might affect their frame properties. (We dealt with some of the effects of discretization in \cite{WilsonACHA}.) Unfortunately, Theorem \ref{wheeden} only partly addresses the problem of errors, because the changes of variable it treats are very special: small (relative to $Q$) translations, combined with close-to-$I_d$ multiparameter dilations centered around $x_Q$. They have no \dq twisting''. 

Our second main result (Theorem \ref{covid} below) extends Theorem \ref{wheeden} to when the $D^{(Q)}$s are replaced by arbitrary close-to-$I_d$ matrices. Everything still rests on the fundamental result in Theorem \ref{onedimbvresult}, but new technical difficulties arise. Briefly, with the diagonal matrices $D^{(Q)}$, the variables all separate, and separate the same way for each $Q$. With arbitrary matrices the variables get mixed up, and mixed up in different ways for different $Q$s. Dealing with this forced us to require a multi-variable bounded variation condition stronger than that of Definition \ref{nbv}.

\begin{definition}\label{snbv}We say that a measurable $f:\Rd\to{\bf C}$ belongs to $SNBV(\Rd)$\footnote{\dq Strong normalized bounded variation on $\Rd$''.} if, for every pair of vectors $(v,y)\in\Rd\times\Rd$, the function $\psi_{v,y}:\R\to{\bf C}$ defined by
$$\psi_{v,y}(t):= f(y+tv)$$
is of bounded variation on $\R$, with $V_{\psi_{v,y}}(\R)\leq 1$. If $Q\subset\Rd$ is a (possibly non-dyadic) cube we say that $f\in SNBV(Q)$ if $f\in SNBV(\Rd)$ and also satisfies $\{x:\ f(x)\not=0\}\subset Q$. If we also have $\int f\,dx=0$ we say $f\in SNBV_0(Q)$. If $\Gee$ is a set of (possibly non-dyadic) cubes in $\Rd$ and $\{f^{(Q)}\}_{Q\in\Gee}$ is a family of functions such that every $f^{(Q)}\in SNBV(Q)$ ($SNBV_0(Q)$) we call $\{f^{(Q)}\}_{Q\in\Gee}$ an $SNBV$ ($SNBV_0$) family.
\end{definition}

\emph{Remark.} $SNBV$ equals $NBV$ plus preservation under affine transformations: $f\in SNBV\iff g(x):=f(b+Ax)\in NBV$ for all $b\in\Rd$ and all $d\times d$ real matrices $A$.\medskip

\emph{Remark.} If $Q\in\Dee_d$ and $L\subset\Rd$ is a line, $L$ meets at most $d+1$ of $Q$'s dyadic children (exercise); therefore every $h^{(Q),k}$ equals $2(d+1)$ times an $SNBV_0(Q)$ function.\medskip

Before stating Theorem \ref{covid} we will say precisely how we are measuring the sizes of matrices and vectors.

\begin{definition}\label{matrixnorm}If $B:=[b_{ij}]$ is a real $d\times l$ matrix then
$$\nrm{B}_\infty:=\max_{1\leq i\leq d}\sum_{j=1}^l \modd{b_{ij}}.$$\end{definition}

\rmk. If $A$ is a square matrix, we will write $\modd A$ for the determinant of $A$.\medskip

\rmk. If $x=(x_1,\ldots,x_d)$, written as a column vector, then $\nrm{x}_\infty=\max_i\modd{x_i}$. If $B$ is a $d\times l$ matrix then $\nrm{B}_\infty$ is $B$'s $(\R^l,\nrm{\cdot}_\infty)\to(\Rd,\nrm{\cdot}_\infty)$ operator norm.

\begin{theorem}\label{covid}There are positive constants $\tilde c(d)$ and $C(d)$, depending only on $d$, such that the following holds. Let $0\leq\eta\leq\tilde c(d)$ and let $\{f^{(Q)}\}_{Q\in\Dee_d}$ be any  $SNBV$ family.  Suppose that, for every $Q\in\Dee_d$, there is a real $d\times d$ matrix $A^{(Q)}$, and there is a vector $y^{(Q)}\in\Rd$, such that
\begin{equation}\label{matrixvectorbound}\nrm{A^{(Q)}-I_d}_\infty+\nrm{y^{(Q)}}_\infty\leq\eta.\end{equation}
For every $Q\in\Dee_d$ define
\begin{equation}\label{originalperturb}\widetilde{f^{(Q)}}(x):=f^{(Q)}(x_Q+A^{(Q)}(\ell(Q)y^{(Q)}+x-x_Q)).\end{equation}
Then
$$\bignrm{\bigbrace{{f^{(Q)}-\modd{A^{(Q)}}\widetilde{f^{(Q)}}\over\modd{Q}^{1/2}}}_{Q\in\Dee_d}}_{AO(\Dee_d)}\leq C(d)\eta^{1/2}.$$
\end{theorem}

Our chief interest in Theorem \ref{covid} lies in the following corollary, which directly applies to perturbations of Haar functions. 

\begin{corollary}\label{covidcor}There are positive constants $\tilde c(d)$ and $C(d)$, depending only on $d$, such that the following holds. Let $0\leq\eta\leq\tilde c(d)$ and let $\{f^{(Q)}\}_{Q\in\Dee_d}$ be any  $SNBV_0$ family.  Suppose that, for every $Q\in\Dee_d$, there is a real $d\times d$ matrix $A^{(Q)}$, and there is a vector $y^{(Q)}\in\Rd$, such that
\begin{equation}\label{matrixvectorbound}\nrm{A^{(Q)}-I_d}_\infty+\nrm{y^{(Q)}}_\infty\leq\eta.\end{equation}
For every $Q\in\Dee_d$ define
\begin{equation}\label{originalperturb}\widetilde{f^{(Q)}}(x):=f^{(Q)}(x_Q+A^{(Q)}(\ell(Q)y^{(Q)}+x-x_Q)).\end{equation}
Then
$$\bignrm{\bigbrace{{f^{(Q)}-\widetilde{f^{(Q)}}\over\modd{Q}^{1/2}}}_{Q\in\Dee_d}}_{AO(\Dee_d)}\leq C(d)\eta^{1/2}.$$
\end{corollary}

\rmk. It is important to note the difference between the theorem and the corollary. The theorem does not ask that the $f^{(Q)}$s have integral 0. \medskip

Corollary \ref{covidcor} immediately implies the following results for the $d$-dimensional Haar functions.

\begin{corollary}\label{covidcor2}Let $\Gamma:=\{(Q,k):\ Q\in\Dee_d,\ 1\leq k<2^d\}$. There are positive constants $\tilde c(d)$ and $C(d)$, depending only on $d$, such that the following holds. Let $0\leq\eta\leq\tilde c(d)$ and let $\{h_{(Q),k}\}_{(Q,k)\in\Gamma}$ be the Haar functions.  Suppose that, for every $(Q,k)\in\Gamma$, there is a real $d\times d$ matrix $A^{(Q),k}$, and there is a vector $y^{(Q),k}\in\Rd$, such that
$$\nrm{A^{(Q),k}-I_d}_\infty+\nrm{y^{(Q),k}}_\infty\leq\eta.$$
For every $(Q,k)\in\Gamma$ define
$$\widetilde{h_{(Q),k}}(x):=h_{(Q),k}(x_Q+A^{(Q),k}(\ell(Q)y^{(Q),k}+x-x_Q)).$$
Then
$$\bignrm{\bigbrace{h_{(Q),k}-\widetilde{h_{(Q),k}}}_{(Q,k)\in\Gamma}}_{AO(\Gamma)}\leq C(d)\eta^{1/2}.$$
\end{corollary}

\begin{corollary}\label{covidcor3}Let $\Gamma:=\{(Q,k):\ Q\in\Dee_d,\ 1\leq k<2^d\}$. There are positive constants $\tilde c(d)$ and $C(d)$, depending only on $d$, such that the following holds. Let $0\leq\eta\leq\tilde c(d)$ and let $\{h_{(Q),k}\}_{(Q,k)\in\Gamma}$ be the Haar functions. Suppose that, for every $(Q,k)\in\Gamma$, there are two real $d\times d$ matrices $A^{(Q),k}$ and $A^{(Q),k}_*$, and there are vectors $y^{(Q),k}$ and $y^{(Q),k}_*$ in $\Rd$, such that
\begin{eqnarray*}\nrm{A^{(Q),k}-I_d}_\infty+\nrm{y^{(Q),k}}_\infty&\leq&\eta\\
\nrm{A^{(Q),k}_*-I_d}_\infty+\nrm{y^{(Q),k}_*}_\infty&\leq&\eta.\end{eqnarray*}
For every $(Q,k)\in\Gamma$ define
\begin{eqnarray*}\widetilde{h_{(Q),k}}(x)&:=&h_{(Q),k}(x_Q+A^{(Q),k}(\ell(Q)y^{(Q),k}+x-x_Q))\\
{h_{(Q),k}^*}(x)&:=&h_{(Q),k}(x_Q+A^{(Q),k}_*(\ell(Q)y^{(Q),k}_*+x-x_Q)).\end{eqnarray*}
Then, for any $f\in L^2$, the series
$$\sum_{(Q,k)\in\Gamma}\langle f,\widetilde{h_{(Q,k)}}\rangle h_{(Q,k)}^*$$
converges unconditionally to some $f^*\in L^2$, and this $f^*$ satisfies $\nrm{f-f^*}_2\leq C(d)\eta^{1/2}\nrm{f}_2$.\end{corollary}

The proof of Theorem \ref{covid} is based on the one for Theorem \ref{wheeden} in \cite{WilsonBvconvao}. The chief new elements are a quantitative form of the $LU$ matrix factorization (see below) and more careful manipulations. 

We give Theorem \ref{genofabg}'s short proof in Section 2. We state and prove lemmas in Section 3. In Section 4 we prove Theorem \ref{covid} and Corollary \ref{covidcor}, and, after some discussion, we use Theorem \ref{covid} to prove the following sharp estimate (Corollary \ref{covidcordyadav}) on the sensitivity of $L^2$ dyadic averages to close-to-the-identity affine changes of variable. We use the familiar notation that, if a measurable $E\subset\Rd$ is bounded and has positive measure, and $g\in\Lwunloc$, then
$$g_E:={1\over\modd E}\int_E g\,dt,$$
the average of $g$ over $E$.

\begin{corollary}\label{covidcordyadav}There is a positive constant $C(d)$ such that the following is true. Let $\{y^{(Q)}\}_{Q\in\Dee_d}$ and $\{A^{(Q)}\}_{Q\in\Dee_d}$ satisfy the hypotheses of Theorem \ref{covid}. For $Q\in\Dee_d$ define $Q^*$ by $\chi_{Q^*}(x):=\chi_Q(x_Q+A^{(Q)}(\ell(Q)y^{(Q)}+x-x_Q))$. Then, for every $g\in L^2(\Rd)$,
$$\bignrm{\left(\sum_{Q\in\Dee_d}\bigmodd{g_Q-g_{Q^*}}^2\chi_Q\right)^{1/2}}_2\leq C(d)\eta^{1/2}\nrm g_2.$$
The exponent $1/2$ on $\eta$ can't be improved.\end{corollary}

We warn the reader that Section 4 is the paper's longest and most technical section.

We sometimes use \lq$C$' to mean a constant that might change from occurrence to occurrence. We will not always state the parameters that $C$ depends on. If $E$ and $F$ are sets, we write $E\subset F$ to mean $E\subseteq F$. 

We indicate the end of a proof with \cs.

\section{Proof of Theorem \ref{genofabg}}

As we observed above, every $h^{(Q),k}$ equals 4 times an $NBV_0(Q)$ function. Therefore, except for that factor of 4,
$$\{h^{(Q),k}\}_{Q\in\Dee_d\atop 1\leq k<2^d}$$
equals a union of $2^d-1$ families like those considered by Theorem \ref{wheeden}. Let us temporarily fix $1\leq k<2^d$. Suppose that, for every $Q\in\Dee_d$, we have a Borel probability measure $\mu^{(Q)}$ on $[-1,1]^d$. For each $Q\in\Dee_d$, define
$$g^{(Q)}(x):=\int_{[-1,1]^d}\left(h^{(Q),k}(x)-h^{(Q),k}(x-\ell(Q)\eta y)\right)\,d\mu^{(Q)}(y),$$
where $0\leq \eta<1/2$. It's clear that each $g^{(Q)}$ has support contained inside $(1+2\eta)Q$. It's also clear that, if the $\nrm\cdot_\infty$ distance between $x$ and any of $h^{(Q),k}$'s discontinuities is $>\eta\ell(Q)$ then 
$$h^{(Q),k}(x)-h^{(Q),k}(x-\ell(Q)\eta y)=0$$
for all $y\in[-1,1]^d$, and $g^{(Q)}(x)=0$. We claim that 
$$\bigbrace{{g^{(Q)}\over\modd{Q}^{1/2}}}_{Q\in\Dee_d}$$
is almost-orthogonal, with an almost-orthogonality constant $\leq C(d)\eta^{1/2}$. Here is the proof. Let $\Eff\subset\Dee_d$ be finite, $\{\lambda_Q\}_{Q\in\Eff}\subset\C$, and $\sum_\Eff\modd{\lambda_Q}^2\leq1$. By Theorem \ref{wheeden}, for every choice of $y:=\{y^{(Q)}\}_{Q\in\Eff}\subset[-1,1]^d$, the finite family
$$\bigbrace{{h^{(Q),k}(x)-h^{(Q),k}(x-\ell(Q)\eta y^{(Q)})\over\modd{Q}^{1/2}}}_{Q\in\Eff}$$
is almost-orthogonal, with an $AO(\Dee_d)$ norm $\leq C(d)\eta^{1/2}$. Let $([-1,1]^d)^\Eff$ denote the (finite) cartesian product $[-1,1]^d\times[-1,1]^d\times\cdots\times[-1,1]^d$ (one factor for each $Q\in\Eff$), and give it the natural product probability measure induced by the $\mu^{(Q)}$s for $Q\in\Eff$. Call this measure $M$. For every $x\in\Rd$,
$$\sum_{Q\in\Eff}\lambda_Q {g^{(Q)}(x)\over\modd{Q}^{1/2}}=\int_{([-1,1]^d)^\Eff}\left(\sum_{Q\in\Eff}\lambda_Q {h^{(Q),k}(x)-h^{(Q),k}(x-\ell(Q)\eta y^{(Q)})\over\modd{Q}^{1/2}}\right)\,dM(y).$$
By the integral form of Minkowski's inequality,
\begin{eqnarray*}&\left(\int_\Rd\bigmodd{\sum_{Q\in\Eff}\lambda_Q {g^{(Q)}(x)\over\modd{Q}^{1/2}}}^2\,dx\right)\\
&\leq\int_{([-1,1]^d)^\Eff}\left(\int_\Rd \bigmodd{\sum_{Q\in\Eff}\lambda_Q {h^{(Q),k}(x)-h^{(Q),k}(x-\ell(Q)\eta y^{(Q)})\over\modd{Q}^{1/2}}}^2\,dx\right)^{1/2}\,dM(y)\\
&\leq\int_{([-1,1]^d)^\Eff} C(d)\eta^{1/2}\,dM(y)=\\
&=C(d)\eta^{1/2},\end{eqnarray*}
proving the claim.

We get Theorem \ref{genofabg}  by taking some non-negative $\psi\in L^1(\Rd)$ with support contained in $[-1,1]^d$ and satisfying $\int\psi\,dx=1$. We set $d\mu^{(Q)}(y'):=\psi(y')\,dy'$. It's then obvious that
\begin{eqnarray*}g^{(Q)}(x)&=&h^{(Q),k}(x)-\int_{[-1,1]^d}h^{(Q),k}(x-\ell(Q)\eta y)\,\psi(y)\,dy\\
&=&h^{(Q),k}(x)-\int_\Rd h^{(Q),k}(x-\ell(Q)\eta y)\,\psi(y)\,dy\\
&=&h^{(Q),k}(x)-\int_\Rd h^{(Q),k}(x-t)\,\eta^{-d}\ell(Q)^{-d}\psi(t/(\eta\ell(Q)))\,dt\\
&=&h^{(Q),k}(x)-(h^{(Q),k}*\psi_{\eta\ell(Q)})(x).\end{eqnarray*}
We do this for every $1\leq k<2^d$ and, if $\eta$ is small enough, get, for any $\epsilon>0$, a family
$$\bigbrace{{h^{(Q),k}*\psi_{\eta\ell(Q)}\over\modd Q^{1/2}}}_{Q\in\Dee_d\atop 1\leq k<2^d}$$
with supports as close to the dyadic cubes as we like (contained in $(1+\epsilon)Q$), equalling the respective Haar functions for $x$s more than $\epsilon\ell(Q)$ distant from their discontinuities, and satisfying
$$(1-C(d)\sqrt\epsilon)^2\nrm f_2^2\leq \sum_{Q\in\Dee_d\atop 1\leq k<2^d}\bigmodd{\left\langle f,{h^{(Q),k}*\psi_{\eta\ell(Q)}\over\modd Q^{1/2}}\right\rangle}^2\leq (1+C(d)\sqrt\epsilon)^2\nrm f_2^2$$
for all $f\in L^2$. It's also clear that this family has the required \dq generating property'' (\ref{genfamily}). We have proved Theorem \ref{genofabg} and, as promised,  a little more. \cs\medskip

\rmk. We note in passing that the proof of Theorem \ref{genofabg} does not use the full power of Theorem \ref{wheeden}. It only requires stability with respect to translations. Dilations play no role.

\section{Preliminary lemmas}

The quantitative form of the $LU$ factorization we mentioned comes from \cite{DM} (Theorem 4.2 (4)). We state it in the form we will use.

\begin{lemma}\label{linalglemma} Let $A$ be a $d\times d$ real matrix and $0\leq\eta<1$. If $\nrm{A-I_d}_\infty\leq\eta$ then we can write $A=LU$, where $L$ is lower triangular, $U$ is upper triangular, and 
$$\max(\nrm{L-I_d}_\infty,\nrm{U-I_d}_\infty)\leq {\eta\over 1-\eta}.$$
In particular, if $\eta\leq1/2$, we have 
$$\max(\nrm{L-I_d}_\infty,\nrm{U-I_d}_\infty)\leq 2\eta.$$
\end{lemma}

\emph{Remark.} In proving Theorem \ref{covid} we will use Lemma \ref{linalglemma} to write the linear transformation defined by $A^{(Q)}$ as a telescoping sum of one-dimensional shear transformations, to which we can apply the methods of \cite{WilsonBvconvao}.\medskip

The preceding lemma implies that, if $\nrm{A-I_d}_\infty\leq\eta\leq1/2$, then $\modd{1-\modd A}\leq \exp(4d\eta)-1$, which is a crude estimate, but good enough for us. 

We recall two elementary lemmas from \cite{WilsonBvconvao}.

\begin{lemma}\label{elementary}(\cite{WilsonBvconvao}, Lemma 1) Let $I\subset \R$ be an interval, $f:I\to{\bf C}$, with $f$ of bounded variation on $I$ with total variation $=:V_f(I)$. Suppose that $b\in L^1(I)$, $b$ is real-valued, and $\int b\,dx=0$. Then
$$\bigmodd{\int_I f(x)\,b(x)\,dx}\leq (1/2)V_f(I)\nrm{b}_1.$$\end{lemma}

\begin{lemma}\label{next}(\cite{WilsonBvconvao}, Lemma 5) If $(a,b)$ and $(a',b')$ are two bounded intervals then
$$\int_\R \modd{\chi_{(a,b)}(x)-\chi_{(a',b')}(x)}\,dx\leq \modd{a-a'}+\modd{b-b'}.$$
\end{lemma}

\emph{Remark.} The same estimate obviously holds if one or both intervals are closed or half-open.\medskip

The next two lemmas require a technical observation. Let $0<\alpha<\infty$ and $\beta\in\R$, and define $\psi(x):= \alpha x+\beta$. Then
$$\chi_{(a,b)}(\psi(x))=\chi_{(\psi^{-1}(a),\psi^{-1}(b))}(x),$$
where $\psi^{-1}(x):=(x-\beta)/\alpha$ is $\psi$'s inverse function. Obviously, we also get
$$\chi_{(a,b)}(\psi^{-1}(x))=\chi_{(\psi(a),\psi(b))}(x).$$

\begin{lemma}\label{diffbound}Let $f:\R\to{\bf C}$ be of bounded variation on \R\ and let $\psi$ be as in the preceding paragraph. Define
$$g(x):= f(x)-\alpha f(\psi(x)).$$
For $J$ a bounded interval define
$$\phi(x):= \chi_{J_l}(x)-\chi_{J_r}(x),$$
where $J_l$ and $J_r$ are $J$'s respective left and right halves. Then:
\begin{eqnarray*}&\bigmodd{\int_\R g(x)\,\phi(x)\,dx}\\ &\leq(1/2)V_f(K(J\cup\psi[J]))\int_\R(\modd{\chi_{J_l}(x)-\chi_{\psi[J_l]}(x)}+\modd{\chi_{J_r}(x)-\chi_{\psi[J_r]}(x)})\,dx,\end{eqnarray*}
where $\psi[S]$ means the the image of a set $S$ under $\psi$, and $K(J\cup\psi[J])$ is the convex hull of $J\cup\psi[J]$.
\end{lemma}

{\bf Proof of Lemma \ref{diffbound}.} After we do a change of variable,
\begin{eqnarray*}\int_\R g(x)\,\phi(x)\,dx&=&\int_\R f(x)\,\phi(x)\,dx-\int_\R f(x)\,\phi(\psi^{-1}(x))\,dx\\
&=&\int_\R f(x)\,(\phi(x)-\phi(\psi^{-1}(x)))\,dx,\end{eqnarray*}
where, because of $\phi$'s and $\psi$'s simple forms, 
$$\int_\R \phi(x)\,dx=\int_\R\phi(\psi^{-1}(x))\,dx=0.$$
We have
$$\phi(x)-\phi(\psi^{-1}(x))=\chi_{J_l}(x)-\chi_{\psi[J_l]}(x) - \left( \chi_{J_r}(x)-\chi_{\psi[J_r]}(x)\right),$$
which is clearly 0 outside $K(J\cup\psi[J])$. Also,
$$\int_\R\modd{\phi(x)-\phi(\psi^{-1}(x))}\,dx\leq \int_\R(\modd{\chi_{J_l}(x)-\chi_{\psi[J_l]}(x)}+\modd{\chi_{J_r}(x)-\chi_{\psi[J_r]}(x)})\,dx.$$
Our result now follows from Lemma \ref{elementary}, with the present $f$ playing the role of \dq$f$'' and $\phi(x)-\phi(\psi^{-1}(x))$ acting as \dq$b$''. \cs

\begin{lemma}\label{excessive}Let $I\subset\R$ be a bounded interval with center $x_I$. Let $\alpha$ and $\tau$ be real numbers such that $\modd{1-\alpha}+\modd\tau \leq\eta$, where $0\leq\eta\leq 1/2$, and set 
$$\psi(x):= \alpha(x-x_I+\tau\ell(I))+x_I.$$
Let $f:\R\to{\bf C}$ be of bounded variation on \R. If $J\subset I$ is any subinterval of $I$ such that $\ell(J)\geq\eta\ell(I)$, then
$$\bigmodd{\int_\R (f(x)-\alpha f(\psi(x)))\,(\chi_{J_l}(x)-\chi_{J_r}(x))\,dx}\leq 2\eta\ell(I) V_f(\widetilde J),$$
where $\widetilde J$ is the concentric triple of $J$, and $J_l$ and $J_r$ mean what they did in Lemma \ref{diffbound}.\end{lemma}

\emph{Remark.} If $\alpha=1$ and $\tau=0$ then $\psi(x)=x$. If $\alpha$ is close to 1 and $\tau$ is close to 0 then $\psi(x)$ is close to $x$ on $I$ and (in a certain sense) $f(x)$ is close to $\alpha f(\psi(x))$. We think of $\alpha f(\psi(x))$ as a perturbation of $f(x)$.\medskip

{\bf Proof of Lemma \ref{excessive}.} If $x\in I$ then
\begin{eqnarray*}\modd{x-\psi(x)}&=&\modd{(\alpha-1)(x-x_I)+\alpha\tau\ell(I)}\\
&\leq&\left({s\over 2}+(1+s)(\eta-s)\right)\ell(I)
\end{eqnarray*}
for some $0\leq s\leq\eta$, where $s:=\modd{\alpha-1}$. If we set $\phi(s):=\left({s\over 2}+(1+s)(\eta-s)\right)=\eta+s\eta-s/2-s^2$, then $\phi'(s)=\eta-1/2-2s\leq0$ for $s\geq 0$ (since $\eta\leq1/2$), so $\phi(s)\leq\phi(0)=\eta$, implying
\begin{equation}\label{psiestimate}\modd{x-\psi(x)}\leq\eta\ell(I).\end{equation}
Therefore, if $J\subset I$ and $\ell(J)\geq\eta\ell(I)$, 
$$K(J\cup\psi[J])\subset \widetilde J,$$
where $K(J\cup\psi[J])$ is as in Lemma \ref{diffbound} above. The result now follows from Lemma \ref{diffbound}, since 
$$\int_\R(\modd{\chi_{J_l}(x)-\chi_{\psi[J_l]}(x)}+\modd{\chi_{J_r}(x)-\chi_{\psi[J_r]}(x)})\,dx\leq
4\eta\ell(I)$$
by Lemma \ref{next} and inequality (\ref{psiestimate}). \cs\medskip

The next lemma describes what happens when a whole family of functions $\{f^{(I)}\}_{I\in\Dee_1}$ like the one in Lemma \ref{excessive} is subjected to a family of perturbations like the one seen there. Its argument, like that of Lemma 2 in \cite{GZ}, relies on inner product estimates and the Schur test. 

\begin{lemma}\label{perturbdiff}Let $0\leq \eta\leq 1/2$. Suppose that, for every $I\in\Dee_1$, we have:\par
a) real numbers $\alpha_I$ and $\tau_I$ such that $\modd{1-\alpha_I}+\modd{\tau_I}\leq\eta$;\par
b) the function $\psi^{(I)}:\R\to\R$, defined by 
$$\psi^{(I)}(x):= \alpha_I(x-x_I+\tau_I\ell(I))+x_I;$$\par
c) a function $f^{(I)}:\R\to{\bf C}$ of bounded variation such that $V_{f^{(I)}}(\R)\leq1$, $\{x:\ f^{(I)}(x)\not=0\}\subset I$, and $\{x:\ f^{(I)}(\psi^{(I)}(x))\not=0\}\subset I$.\par
Define, for each $I\in\Dee_1$,
\begin{equation}\label{diffofgs}g^{(I)}(x):= f^{(I)}(x)-\alpha_I f^{(I)}(\psi^{(I)}(x)).\end{equation}
Then 
\begin{equation}\label{alterbessel}\bigbrace{{g^{(I)}(x)\over\modd{I}^{1/2}}}_{\Dee_1}\end{equation}
is a Bessel family in $L^2(\R)$ with a Bessel bound $\leq C\eta$, where $C$ is an absolute constant.\end{lemma}

\rmk. Alternatively, we could say that (\ref{alterbessel}) is almost-orthogonal, with an almost-orthogonality constant (or norm) $\leq C\eta^{1/2}$.\medskip

\emph{Remark.} Notice that Lemma \ref{perturbdiff} does {not} assume $\int f^{(I)}\,dx=0$. The necessary cancelation comes by taking differences (\ref{diffofgs}). This will be important later.\medskip

{\bf Proof of Lemma \ref{perturbdiff}.} We will apply the Schur test after we get estimates for the $L^2$ inner products
$$\langle g^{(I)},h^{(J)}\rangle$$
for arbitrary dyadic intervals $I$ and $J$, where the functions $h^{(J)}$ are as defined in (\ref{haarfunctiondef}).

We observe that $\{x:\ g^{(I)}(x)\not=0\}\subset I$ and $\int g^{(I)}(x)\,dx=0$. Therefore
$$\langle g^{(I)},h^{(J)}\rangle=0$$ 
unless $J\subset I$. For $J\subset I$ we consider two cases: $\ell(J)\geq \eta\ell(I)$ and $\ell(J)<\eta\ell(I)$. 

In the first case, by Lemma \ref{excessive}, 
$$\bigmodd{\langle g^{(I)},h^{(J)}\rangle}\leq 2\eta\ell(I)V_{f^{(I)}}(\widetilde J),$$
which implies that
$$\bigmodd{\left\langle {g^{(I)}\over\modd{I}^{1/2}},{h^{(J)}\over\modd{J}^{1/2}}\right\rangle}\leq 2\eta\sqrt{{\modd I\over\modd J}}V_{f^{(I)}}(\widetilde J).$$

In the second case, by Lemma \ref{elementary},
\begin{eqnarray*}\bigmodd{\langle g^{(I)},h^{(J)}\rangle}&\leq& (1/2)V_{g^{(I)}}(J)\nrm{h^{(J)}}_1\\
&=&(1/2)V_{g^{(I)}}(J)\modd J,\end{eqnarray*}
implying
$$\bigmodd{\left\langle {g^{(I)}\over\modd{I}^{1/2}},{h^{(J)}\over\modd{J}^{1/2}}\right\rangle}\leq (1/2)V_{g^{(I)}}(J)\sqrt{{\modd J\over\modd I}}.$$

We are now ready to apply the Schur test. First we fix $I\in\Dee_1$ and consider
\begin{eqnarray*}&\sum_{J\in\Dee_1}\bigmodd{\left\langle {g^{(I)}\over\modd{I}^{1/2}},{h^{(J)}\over\modd{J}^{1/2}}\right\rangle}\\
&=\sum_{J\in\Dee_1 \atop J\subset I, \ell(J)\geq\eta\ell(I)}\bigmodd{\left\langle {g^{(I)}\over\modd{I}^{1/2}},{h^{(J)}\over\modd{J}^{1/2}}\right\rangle}+
\sum_{J\in\Dee_1\atop J\subset I, \ell(J)<\eta\ell(I)}\bigmodd{\left\langle {g^{(I)}\over\modd{I}^{1/2}},{h^{(J)}\over\modd{J}^{1/2}}\right\rangle}\\
&=:(I)+(II).\end{eqnarray*}
Let $N$ be the unique natural number such that 
$$2^{-N}\geq\eta>2^{-N-1}.$$

Then
\begin{eqnarray*}(I)&=&\sum_{J\in\Dee_1,\,J\subset I\atop \ell(J)\geq2^{-N}\ell(I)}\bigmodd{\left\langle {g^{(I)}\over\modd{I}^{1/2}},{h^{(J)}\over\modd{J}^{1/2}}\right\rangle}\\
&=&\sum_{k=0}^N\sum_{J\in\Dee_1,\,J\subset I\atop \ell(J)=2^{-k}\ell(I)}\bigmodd{\left\langle {g^{(I)}\over\modd{I}^{1/2}},{h^{(J)}\over\modd{J}^{1/2}}\right\rangle}\\
&\leq&\sum_{k=0}^N\sum_{J\in\Dee_1,\,J\subset I\atop \ell(J)=2^{-k}\ell(I)} 2\eta\sqrt{{\modd I\over\modd J}}V_{f^{(I)}}(\widetilde J)\\
&=&2\eta\sum_{k=0}^N 2^{k/2}\sum_{J\in\Dee_1,\,J\subset I\atop \ell(J)=2^{-k}\ell(I)}V_{f^{(I)}}(\widetilde J).\end{eqnarray*}
But, for each $k$,
$$\sum_{J\in\Dee_1,\,J\subset I\atop \ell(J)=2^{-k}\ell(I)}V_{f^{(I)}}(\widetilde J)\leq
3V_{f^{(I)}}(\R)\leq 3,$$
implying
$$(I)\leq 6\eta\sum_{k=0}^N 2^{k/2}\leq C\eta 2^{N/2}\leq C\eta^{1/2},$$
from the definition of $N$. 

On the other hand,
\begin{eqnarray*}(II)&=&\sum_{J\in\Dee_1,\,J\subset I\atop \ell(J)<2^{-N}\ell(I)}\bigmodd{\left\langle {g^{(I)}\over\modd{I}^{1/2}},{h^{(J)}\over\modd{J}^{1/2}}\right\rangle}\\
&=&\sum_{k=N+1}^\infty\sum_{J\in\Dee_1,\,J\subset I\atop \ell(J)=2^{-k}\ell(I)}\bigmodd{\left\langle {g^{(I)}\over\modd{I}^{1/2}},{h^{(J)}\over\modd{J}^{1/2}}\right\rangle}\\
&\leq&\sum_{k=N+1}^\infty\sum_{J\in\Dee_1,\,J\subset I\atop \ell(J)=2^{-k}\ell(I)}(1/2)V_{g^{(I)}}(J)\sqrt{{\modd J\over\modd I}} \\
&=&(1/2)\sum_{k=N+1}^\infty 2^{-k/2}\sum_{J\in\Dee_1,\,J\subset I\atop \ell(J)=2^{-k}\ell(I)}V_{g^{(I)}}(J).\end{eqnarray*}
For each $k$,
$$\sum_{J\in\Dee_1,\,J\subset I\atop \ell(J)=2^{-k}\ell(I)}V_{g^{(I)}}(J)\leq V_{g^{(I)}}(\R).$$
But
\begin{equation}\label{vgrbound}V_{g^{(I)}}(\R)\leq (1+\modd{\alpha_I})V_{f^{(I)}}(\R)\leq 3.\end{equation}
Therefore
\begin{eqnarray*}(II)&\leq&(3/2)\sum_{k=N+1}^\infty 2^{-k/2}\\
&\leq& C2^{-N/2}\\
&\leq& C\eta^{1/2},\end{eqnarray*}
yielding
\begin{equation}\label{schurbnd1}(I)+(II)\leq C\eta^{1/2}.\end{equation}

Now we go the other way. Fix $J\in\Dee_1$ and consider
\begin{eqnarray*}&\sum_{I\in\Dee_1}\bigmodd{\left\langle {g^{(I)}\over\modd{I}^{1/2}},{h^{(J)}\over\modd{J}^{1/2}}\right\rangle}\\
&=\sum_{I\in\Dee_1 \atop I\supset J, \ell(J)\geq\eta\ell(I)}\bigmodd{\left\langle {g^{(I)}\over\modd{I}^{1/2}},{h^{(J)}\over\modd{J}^{1/2}}\right\rangle}+
\sum_{I\in\Dee_1\atop I\supset J, \ell(J)<\eta\ell(I)}\bigmodd{\left\langle {g^{(I)}\over\modd{I}^{1/2}},{h^{(J)}\over\modd{J}^{1/2}}\right\rangle}\\
&=:(I')+(II').\end{eqnarray*}
Then, as before,
\begin{eqnarray*}(I')&=&\sum_{I\in\Dee_1 \atop I\supset J, \ell(J)\geq\eta\ell(I)}\bigmodd{\left\langle {g^{(I)}\over\modd{I}^{1/2}},{h^{(J)}\over\modd{J}^{1/2}}\right\rangle}\\
&=&\sum_{I\in\Dee_1 \atop I\supset J, \ell(I)\leq\eta^{-1}\ell(J)}\bigmodd{\left\langle {g^{(I)}\over\modd{I}^{1/2}},{h^{(J)}\over\modd{J}^{1/2}}\right\rangle}\\
&=&\sum_{I\in\Dee_1 \atop I\supset J, \ell(I)\leq 2^N\ell(J)}\bigmodd{\left\langle {g^{(I)}\over\modd{I}^{1/2}},{h^{(J)}\over\modd{J}^{1/2}}\right\rangle}\\
&=&\sum_{k=0}^N\sum_{I\in\Dee_1,\,I\supset J\atop \ell(I)= 2^k\ell(J)}\bigmodd{\left\langle {g^{(I)}\over\modd{I}^{1/2}},{h^{(J)}\over\modd{J}^{1/2}}\right\rangle}\\
&\leq&\sum_{k=0}^N\sum_{I\in\Dee_1,\,I\supset J\atop \ell(I)=2^k\ell(J)} 2\eta\sqrt{{\modd I\over\modd J}}V_{f^{(I)}}(\widetilde J)\\
&=&2\eta\sum_{k=0}^N 2^{k/2}\sum_{I\in\Dee_1,\,I\supset J\atop \ell(I)=2^k\ell(J)}V_{f^{(I)}}(\widetilde J)\\
&\leq&2\eta\sum_{k=0}^N 2^{k/2}\\
&\leq& C\eta^{1/2},\end{eqnarray*}
where
$$\sum_{I\in\Dee_1,\,I\supset J\atop \ell(I)=2^k\ell(J)}V_{f^{(I)}}(\widetilde J)\leq 1$$
because $V_{f^{(I)}}(\widetilde J)\leq V_{f^{(I)}}(\R)\leq1$ and, for every $k\geq0$, every $J\in\Dee_1$ is contained in only one $I\in\Dee_1$ with $\ell(I)=2^k\ell(J)$. That takes care of $(I')$. 

For $(II')$, we again fix $J\in\Dee_1$, and write
\begin{eqnarray*}(II')&=&\sum_{I\in\Dee_1\atop I\supset J, \ell(J)<\eta\ell(I)}\bigmodd{\left\langle {g^{(I)}\over\modd{I}^{1/2}},{h^{(J)}\over\modd{J}^{1/2}}\right\rangle}\\
&=&\sum_{I\in\Dee_1\atop I\supset J, \ell(I)>\eta^{-1}\ell(J)}\bigmodd{\left\langle {g^{(I)}\over\modd{I}^{1/2}},{h^{(J)}\over\modd{J}^{1/2}}\right\rangle}\\
&=&\sum_{I\in\Dee_1\atop I\supset J, \ell(I)>2^N\ell(J)}\bigmodd{\left\langle {g^{(I)}\over\modd{I}^{1/2}},{h^{(J)}\over\modd{J}^{1/2}}\right\rangle}\\
&=&\sum_{k=N+1}^\infty\sum_{I\in\Dee_1,\,I\supset J\atop \ell(I)=2^k\ell(J)}\bigmodd{\left\langle {g^{(I)}\over\modd{I}^{1/2}},{h^{(J)}\over\modd{J}^{1/2}}\right\rangle}\\
&\leq&\sum_{k=N+1}^\infty\sum_{I\in\Dee_1,\,I\supset J\atop \ell(I)=2^k\ell(J)}(1/2)V_{g^{(I)}}(J)\sqrt{{\modd J\over\modd I}} \\
&\leq& C\sum_{k=N+1}^\infty 2^{-k/2}\\
&\leq& C2^{-N/2}\\
&\leq& C\eta^{1/2},\end{eqnarray*}
where
$$\sum_{I\in\Dee_1,\,I\supset J\atop \ell(I)=2^k\ell(J)}(1/2)V_{g^{(I)}}(J)\sqrt{{\modd J\over\modd I}}\leq C2^{-k/2}$$
for three reasons:
$$\sqrt{{\modd J\over\modd I}}=2^{-k/2}$$
and
$$V_{g^{(I)}}(J)\leq V_{g^{(I)}}(\R)\leq 3$$
(see (\ref{vgrbound})) and, for every $k\geq0$, every $J\in\Dee_1$ is contained in only one $I\in\Dee_1$ with $\ell(I)=2^k\ell(J)$. We have proved that 
\begin{equation}\label{schurbnd2}(I')+(II')\leq C\eta^{1/2}.\end{equation}

Putting (\ref{schurbnd1}) and (\ref{schurbnd2}) together, we have that, for all $I\in\Dee_1$,
$$\sum_{J\in\Dee_1}\bigmodd{\left\langle {g^{(I)}\over\modd{I}^{1/2}},{h^{(J)}\over\modd{J}^{1/2}}\right\rangle}\leq C\eta^{1/2}$$
and, for all $J\in\Dee_1$,
$$\sum_{I\in\Dee_1}\bigmodd{\left\langle {g^{(I)}\over\modd{I}^{1/2}},{h^{(J)}\over\modd{J}^{1/2}}\right\rangle}\leq C\eta^{1/2}.$$
Consider the mapping $T:\ell^2(\Dee_1)\to\ell^2(\Dee_1)$, initially defined for finite sums, given by
$$T\left(\{\lambda_I\}_{I\in\Dee_1}\right):=\bigbrace{\sum_{I\in\Dee_1}\lambda_I\left\langle {g^{(I)}\over\modd{I}^{1/2}},{h^{(J)}\over\modd{J}^{1/2}}\right\rangle}_{J\in\Dee_1}.$$
By the Schur test, $T$'s operator norm is $\leq C\eta^{1/2}$. Therefore $T$ extends to a bounded operator (with the same norm) on all of $\ell^2(\Dee_1)$. 

Let $\Eff\subset\Dee_1$ be finite, and consider
$$f=\sum_{I\in\Eff}\lambda_I {g^{(I)}\over\modd{I}^{1/2}}.$$
Since the Haar functions form a complete orthonormal system in $L^2(\R)$,
\begin{eqnarray*}\int\modd{f}^2\,dx&=&\sum_{J\in\Dee_1}\bigmodd{\left\langle f,{h^{(J)}\over\modd{J}^{1/2}}\right\rangle}^2\\
&=&\sum_{J\in\Dee_1}\bigmodd{\sum_{I\in\Eff} \lambda_I \left\langle {g^{(I)}\over\modd{I}^{1/2}},{h^{(J)}\over\modd{J}^{1/2}}\right\rangle}^2\\
&=&\nrm{T\left(\{\lambda_I\}_{I\in\Eff}\right)}_{\ell^2(\Dee_1)}^2\\
&\leq& C\eta\sum_{I\in\Eff}\modd{\lambda_I}^2,\end{eqnarray*}
and therefore
$$\bignrm{\sum_{I\in\Eff}\lambda_I {g^{(I)}\over\modd{I}^{1/2}}}_2\leq C\eta^{1/2}\left(\sum_{I\in\Eff}\modd{\lambda_I}^2\right)^{1/2},$$
proving Lemma \ref{perturbdiff}. \cs

\section{Proofs of Theorem \ref{covid} and consequences}

We will first prove Theorem \ref{covid} under the assumption that, for each $Q\in\Dee_d$, $\{x:\ f^{(Q)}(x)\not=0\}$ lies inside the central third of $Q$; i.e., the unique cube whose concentric triple equals $Q$. 

We set $\tilde c(d):=1/(20d)$. 

We claim that, given (\ref{matrixvectorbound}) and our choice of $\tilde c(d)$, every $\widetilde{f^{(Q)}}\in SNBV(Q)$. (We will use this fact near the end of the proof.) Here is the proof of the claim. 
The bounded variation condition is immediate. We need to check the support. We recall that, if $A^{(Q)}$ is a $d\times d$ real matrix and $x\in\Rd$, then $\nrm{A^{(Q)}x}_\infty\leq \nrm{A^{(Q)}}_\infty\nrm{x}_\infty$, and therefore
\begin{eqnarray*}\nrm{A^{(Q)}x}_\infty&=&\nrm{(I_d-(I_d-A^{(Q)}))x}_\infty\\
&\geq& \nrm{x}_\infty-\nrm{I_d-A^{(Q)}}_\infty\nrm{x}_\infty\\
&=&(1-\nrm{I_d-A^{(Q)}}_\infty)\nrm{x}_\infty\\
&\geq& (1-\eta)\nrm{x}_\infty.\end{eqnarray*}
If $y^{(Q)}\in\Rd$, $\nrm{I_d-A^{(Q)}}_\infty+\nrm{y^{(Q)}}_\infty\leq\eta\leq\tilde c(d)$, and  $\nrm{x-x_Q}_\infty\geq (1/2)\ell(Q)$ then
\begin{eqnarray*}\nrm{(A^{(Q)}(x-x_Q+\ell(Q)y^{(Q)})+x_Q)-x_Q}_\infty&\geq& (1-\eta)\nrm{x-x_Q+\ell(Q)y^{(Q)}}_\infty\\
&\geq& (1-\eta)(\ell(Q)/2-\eta\ell(Q))\\
&>&\ell(Q)/3.\end{eqnarray*}
But $x\notin Q\Rightarrow \nrm{x-x_Q}_\infty\geq (1/2)\ell(Q)$, which---for this $\eta$---implies $A^{(Q)}(x-x_Q+\ell(Q)y^{(Q)})+x_Q\notin$ $Q$'s central third, since the $\nrm{\cdot}_\infty$-distance between it and $x_Q$ is bigger than $\ell(Q)/3>\ell(Q)/6$. Therefore, having $x$ outside $Q$ forces $\widetilde{f^{(Q)}}(x)=f^{(Q)}(A^{(Q)}(x-x_Q+\ell(Q)y^{(Q)})+x_Q)=0$, proving the claim. 

Because of the claim, the $d=1$ case of Theorem \ref{covid}---with our special assumption---follows immediately from Lemma \ref{perturbdiff}. Therefore we assume that $d\geq2$.

For future reference we note that, under the hypotheses of Theorem \ref{covid}, if $\nrm{x-x_Q}_\infty\geq\ell(Q)/2$ then 
$$\nrm{A^{(Q)}(x-x_Q)}_\infty\geq \ell(Q)/3.$$
We'll use this near the end of the proof too.

We now make a second preliminary assumption: every $y^{(Q)}=0$. We do this because the hard work of proving Theorem \ref{covid} comes from dealing with the $A^{(Q)}$s, for which the $y^{(Q)}$s are a distraction. We shall see that bringing the $y^{(Q)}$s back in at the end is no trouble. With our two preliminary assumptions, we need to show that
$$\bigbrace{{f^{(Q)}-\modd{A^{(Q)}}\widetilde{f^{(Q)}}\over\modd{Q}^{1/2}}}_{Q\in\Dee_d}$$
is almost-orthogonal with a suitable constant, where 
$$\widetilde{f^{(Q)}}(x):=f^{(Q)}(A^{(Q)}(x-x_Q)+x_Q).$$

We apply the factorization lemma (Lemma \ref{linalglemma}). We factor every $A^{(Q)}$ as $L^{(Q)}U^{(Q)}$, where $L^{(Q)}$ is lower triangular, $U^{(Q)}$ is upper triangular, and both are close to the identity: $\nrm{I_d-L^{(Q)}}_\infty$ and $\nrm{I_d-U^{(Q)}}_\infty$ are $\leq 2\eta$. The crude determinant estimate implies that $\modd{1-\modd{L^{(Q)}}}$ and $\modd{1-\modd{U^{(Q)}}}$ are $\leq Cd\eta\leq C$, an absolute constant, so that $\modd{L^{(Q)}}$ and $\modd{U^{(Q)}}$ are also bounded by absolute constants. Below we will observe that the same applies to certain matrices we will build from $L^{(Q)}$ and $U^{(Q)}$.

Let us define $g^{(Q)}(x):= f^{(Q)}(x_Q+L^{(Q)}x)$. Then, since $y^{(Q)}=0$, we have
\begin{eqnarray*}\modd{A^{(Q)}}\widetilde{f^{(Q)}}(x)&=&\modd{A^{(Q)}}f^{(Q)}(A^{(Q)}(x-x_Q)+x_Q)\\
&=&\modd{A^{(Q)}}f^{(Q)}(L^{(Q)}U^{(Q)}(x-x_Q)+x_Q)\\
&=&\modd{A^{(Q)}}g^{(Q)}(U^{(Q)}(x-x_Q)).\end{eqnarray*}
Therefore
\begin{eqnarray*}\modd{A^{(Q)}}\widetilde{f^{(Q)}}(x)-f^{(Q)}(x)&=&\modd{A^{(Q)}}g^{(Q)}(U^{(Q)}(x-x_Q))-f^{(Q)}(x)\\
&=&\left(\modd{A^{(Q)}}g^{(Q)}(U^{(Q)}(x-x_Q))-\modd{L^{(Q)}}g^{(Q)}(x-x_Q)\right)\\&+&\left(\modd{L^{(Q)}}g^{(Q)}(x-x_Q)-f^{(Q)}(x)\right)\\
&=&\left(\modd{A^{(Q)}}g^{(Q)}(U^{(Q)}(x-x_Q))-\modd{L^{(Q)}}g^{(Q)}(x-x_Q)\right)\\&+&\left(\modd{L^{(Q)}}f^{(Q)}(x_Q+L^{(Q)}(x-x_Q))-f^{(Q)}(x)\right)\\
&=:& I^{(Q)}(x)+II^{(Q)}(x),\end{eqnarray*}
where
\begin{eqnarray*}I^{(Q)}(x)&:=&\modd{A^{(Q)}}g^{(Q)}(U^{(Q)}(x-x_Q))-\modd{L^{(Q)}}g^{(Q)}(x-x_Q)\\
II^{(Q)}(x)&:=&\modd{L^{(Q)}}f^{(Q)}(x_Q+L^{(Q)}(x-x_Q))-f^{(Q)}(x).\end{eqnarray*}

Our job now is to show that the families 
\begin{equation}\label{firstfamily}\bigbrace{{I^{(Q)}\over\modd{Q}^{1/2}}}_{Q\in\Dee_d}\end{equation}
and
\begin{equation}\label{secondfamily}\bigbrace{{II^{(Q)}\over\modd{Q}^{1/2}}}_{Q\in\Dee_d}\end{equation}
are almost-orthogonal, with appropriate constants. 

We remind the reader that 
\begin{eqnarray*}\modd{A^{(Q)}}&=&\modd{L^{(Q)}}\modd{U^{(Q)}}\\
&=&\left(\prod_1^d L^{(Q)}(i,i)\right)\left(\prod_1^d U^{(Q)}(i,i)\right),\end{eqnarray*}
which we will use shortly.

We bound (\ref{firstfamily}) first. Write
\begin{equation}\label{iqsum}I^{(Q)}(x)=\modd{L^{(Q)}}\sum_{k=1}^d \left(\modd{U^{(Q)}_{k-1}}g^{(Q)}(U^{(Q)}_{k-1}(x-x_Q))-\modd{U^{(Q)}_{k}}g^{(Q)}(U^{(Q)}_{k}(x-x_Q))\right),\end{equation}
where
\begin{eqnarray*}U^{(Q)}_{0}&:=&U^{(Q)}\\
U^{(Q)}_{d}&:=&I_d,\end{eqnarray*}
and, for $0<k<d$,
$$U^{(Q)}_k(i,j)=\cases{\delta(i,j)&if $\min(i,j)\leq k$;\cr
U^{(Q)}(i,j)&if $\min(i,j)>k$.\cr}$$
What does this look like? If $M$ is a $d\times d$ matrix then the set of entries $\{M(i,j):\ \min(i,j)\leq k\}$ is the union of $M$'s first $k$ rows (counting from the top) and first $k$ columns (counting from the left), and  $\{M(i,j):\ \min(i,j)>k\}$ is everything else: the entries of a $(d-k)\times (d-k)$ submatrix in $M$'s lower right corner. As $k$ increases, $U^{(Q)}_k$ changes from $U^{(Q)}$ into $I_d$, one elbow-shaped layer at a time, with the elbow joints both opening toward and moving toward the lower right (southeastward). 

Observe that, for every $0\leq k\leq d$, $\nrm{I_d-U^{(Q)}_k}_\infty\leq\nrm{I_d-U^{(Q)}}_\infty$. This is because, on the main diagonal, $U^{(Q)}_k(i,i)=U^{(Q)}(i,i)$ or $U^{(Q)}_k(i,i)=1$; while, off the diagonal, $U^{(Q)}_k(i,j)=U^{(Q)}(i,j)$ or $U^{(Q)}_k(i,j)=0$. We point out (see above) that the determinants of the $U^{(Q)}_k$s, like that of $U^{(Q)}$, are bounded by an absolute constant independent of $k$, $Q$, and $d$.

Suppose that $z=\{z_j\}_1^d\in\Rd$ and $\xi=\{\xi_i\}_1^d:=U^{(Q)}_kz$ (both written as column vectors). Obviously, if $k=0$ then $\xi_i=\sum_{j=1}^d U^{(Q)}(i,j)z_j$, and if $k=d$ then $\xi_i=z_i$ $(1\leq i\leq d$). If $0<k<d$ then
$$\xi_i=\cases{z_i&if $i\leq k$;\cr \sum_{j=i}^d U^{(Q)}(i,j)z_j=\sum_{j=1}^d U^{(Q)}(i,j)z_j&if $i>k$,\cr}$$
where the equality
$$\sum_{j=i}^d U^{(Q)}(i,j)z_j=\sum_{j=1}^d U^{(Q)}(i,j)z_j$$
holds for $i>k$ because $U^{(Q)}$ is upper triangular.

Let's look at one of the terms in (\ref{iqsum}):
\begin{equation}\label{iqterm}\modd{U^{(Q)}_{k-1}}g^{(Q)}(U^{(Q)}_{k-1}(x-x_Q))-\modd{U^{(Q)}_{k}}g^{(Q)}(U^{(Q)}_{k}(x-x_Q)).\end{equation}
Now, 
$$\modd{U^{(Q)}_{k}}=\cases{\prod_{k+1}^d U^{(Q)}(i,i)&if $0\leq k<d$;\cr
1&if $k=d$.\cr}$$
Therefore we can rewrite (\ref{iqterm}) as
\begin{equation}\label{iqtermrewr}\modd{U^{(Q)}_{k}}\left(U^{(Q)}(k,k)g^{(Q)}(U^{(Q)}_{k-1}(x-x_Q))-g^{(Q)}(U^{(Q)}_{k}(x-x_Q))\right).\end{equation}
The determinant at the far left of (\ref{iqtermrewr}) {is bounded by an absolute constant}: we will ignore it. For ease of reading, we write $x-x_Q=z=(z_1,\ldots,z_d)$. Then:
\begin{eqnarray*}&U^{(Q)}(k,k)g^{(Q)}(U^{(Q)}_{k-1}z)=U^{(Q)}(k,k)\times\\
&g^{(Q)}(z_1,z_2,\ldots,z_{k-1},\sum_{j=k}^d U^{(Q)}(k,j)z_j,\sum_{j=k+1}^d U^{(Q)}(k+1,j)z_j,\ldots,U^{(Q)}(d,d)z_d),\end{eqnarray*}
where it's understood that the first $k-1$ entries---$z_1,\,z_2,\,\ldots\,,\,z_{k-1}$---in $g^{(Q)}$'s argument \emph{will not be there} if $k=1$ (because then $k-1=0$). 
It might be easier to see what is going on here if we arrange the entries of $g^{(Q)}$'s argument vertically. We get
\begin{eqnarray*}&z_1\\
&z_2\\
&\vdots\\
&z_{k-1}\\
&\sum_{j=k}^d U^{(Q)}(k,j)z_j\\
&\sum_{j=k+1}^d U^{(Q)}(k+1,j)z_j\\
&\vdots\\
&U^{(Q)}(d,d)z_d\\\end{eqnarray*}
if $k>1$ and 
\begin{eqnarray*}&\sum_{j=1}^d U^{(Q)}(1,j)z_j\\
&\sum_{j=2}^d U^{(Q)}(2,j)z_j\\
&\vdots\\
&U^{(Q)}(d,d)z_d\\\end{eqnarray*}
if $k=1$. In either case, the $k^{th}$ entry in the $U^{(Q)}(k,k)g^{(Q)}(U^{(Q)}_{k-1}z)$ is
$$U^{(Q)}(k,k)z_k+\sum_{j=k+1}^d U^{(Q)}(k,j)z_j$$
(where the sum is empty if $k=d$), and it is (very important!) the only entry containing the variable $z_k$. 

Let's look at the other piece of (\ref{iqtermrewr}) (once again ignoring the determinant way off on the left). It's
\begin{eqnarray*}\label{closeanalyze2}
&g^{(Q)}(U^{(Q)}_{k}z)=\\
&g^{(Q)}(z_1,\ldots,z_k,\sum_{j=k+1}^d U^{(Q)}(k+1,j)z_j,\ldots,U^{(Q)}(d,d)z_d)\end{eqnarray*}
unless $k=d$, in which case it's just $g^{(Q)}(z)$. Once again, we arrange the argument's entries in a column.
\begin{eqnarray*}&z_1\\
&z_2\\
&\vdots\\
&z_{k}\\
&\sum_{j=k+1}^d U^{(Q)}(k+1,j)z_j\\
&\vdots\\
&U^{(Q)}(d,d)z_d,\end{eqnarray*}
unless $k=d$, when we get
\begin{eqnarray*}&z_1\\
&z_2\\
&\vdots\\
&z_d.\end{eqnarray*}

We see that $g^{(Q)}(U^{(Q)}_{k-1}z)$ and 
$g^{(Q)}(U^{(Q)}_{k}z)$ are identical except for their $k^{th}$ entries ($1\leq k\leq d$) and, for both of them, \emph{the variable $z_k$ occurs only there}. This implies that $g^{(Q)}(U^{(Q)}_{k-1}(x-x_Q))$ and $g^{(Q)}(U^{(Q)}_{k}(x-x_Q))$, when written as functions of $x=(x_1,\ldots,x_d)$, are identical
except for their $k^{th}$ entries, and both have the variable $x_k$ occurring only there. 

Before we get too technical it will be useful to note a property shared by all the functions $g^{(Q)}(U^{(Q)}_{j}(x-x_Q))$ ($0\leq j\leq d$): If $x\notin Q$ then $g^{(Q)}(U^{(Q)}_{j}(x-x_Q))=0$. This follows from the inequalities $\nrm{I_d-L^{(Q)}}_\infty\leq2\eta\leq 1/(10d)$ and $\nrm{I_d-U_j^{(Q)}}_\infty\leq2\eta\leq 1/(10d)$ ($0\leq j\leq d$), and our earlier observations. Simply put, $x\notin Q$ implies $\nrm{x-x_Q}_\infty\geq\ell(Q)/2$, implies $\nrm{U^{(Q)}_{k}(x-x_Q)}_\infty\geq(1-2\eta)\ell(Q)/2$, implies $\nrm{L^{(Q)}U^{(Q)}_{j}(x-x_Q)}_\infty\geq(1-2\eta)^2\ell(Q)/2> \ell(Q)/6$ because $2\eta\leq 1/(10d)$, and 
$$g^{(Q)}(U^{(Q)}_{j}(x-x_Q))=f^{(Q)}(x_Q+L^{(Q)}U^{(Q)}_{j}(x-x_Q))=0,$$
because $f^{(Q)}$ is zero outside the central third of $Q$. The nature of the $\nrm{\cdot}_\infty$ norm implies that $g^{(Q)}(U^{(Q)}_{j}(x-x_Q))=0$ if $\modd{x_l-(x_Q)_l}\geq\ell(Q)/2$ for a single $l$.

Suppose we temporarily fix all the variables except $x_k$. Then, as a function of $x_k$, we can write $g^{(Q)}(U^{(Q)}_{k}(x-x_Q))$ as $\phi(x_k)$ for some function $\phi:\R\to{\bf C}$. This function $\phi$ has two interesting properties. It is of bounded variation, with total variation $\leq 1$; and it is zero outside of $I_k(Q)$, the $k^{th}$ cartesian factor of $Q$. The first property is an immediate consequence of the $SNBV$ property. The second was shown in the preceding paragraph. Our observation at the end of that paragraph implies that if $\modd{x_l-(x_Q)_l}\geq\ell(Q)/2$ for some $l\not=k$ then $\phi(x_k)\equiv0$.

The other function we are presently interested in, namely $U^{(Q)}(k,k)g^{(Q)}(U^{(Q)}_{k-1}(x-x_Q))$, when considered as a function of $x_k$ only (recall that the other variables are fixed for now), can be expressed in terms of the same $\phi$; namely, as $\alpha\phi(\alpha(x_k-(x_Q)_k+\tau\ell(Q))+(x_Q)_k)$, where $\alpha=U^{(Q)}(k,k)$ and $\tau$ is a number such that {\bf either} $\modd{\alpha-1}+\modd\tau\leq 7\eta/2$ {\bf or} $\phi(x_k)$ and $\alpha\phi(\alpha(x_k-(x_Q)_k+\tau\ell(Q))+(x_Q)_k)$ are both identically 0. To see these facts, compare the $k^{th}$ entries of $g^{(Q)}(U^{(Q)}_{k}(x-x_Q))$ and $g^{(Q)}(U^{(Q)}_{k-1}(x-x_Q))$. For $g^{(Q)}(U^{(Q)}_{k}(x-x_Q))$ this $k^{th}$ entry is $x_k-(x_Q)_k$. For $g^{(Q)}(U^{(Q)}_{k-1}(x-x_Q))$ it's
\begin{eqnarray}\label{shift}&\alpha(x_k-(x_Q)_k)+\sum_{j=k+1}^d U^{(Q)}(k,j)(x_j-(x_Q)_j)\nonumber\\
&=\alpha\left(x_k-(x_Q)_k+\alpha^{-1}\sum_{j=k+1}^d U^{(Q)}(k,j)(x_j-(x_Q)_j)\right).\end{eqnarray}
We have $\modd{\alpha-1}\leq2\eta$, $0<\alpha^{-1}\leq 1.1$ (since $\eta\leq 1/(20d)$), and 
$$\sum_{j=k+1}^d \modd{U^{(Q)}(k,j)}\leq 2\eta,$$
since $\nrm{I_d-U^{(Q)}}_\infty\leq2\eta$. 
If our functions don't both vanish for all $x_k\in\R$, then no $\modd{x_j-(x_Q)_j}$ ($j\not=k$) can be $\geq\ell(Q)/2$, forcing
$$\bigmodd{\alpha^{-1}\sum_{j=k+1}^d U^{(Q)}(k,j)(x_j-(x_Q)_j)}\leq 2.2\eta\ell(Q)/2=1.1\eta\ell(Q),$$
so we can write
$$\alpha^{-1}\sum_{j=k+1}^d U^{(Q)}(k,j)(x_j-(x_Q)_j)=:\tau\ell(Q)$$
for some $\tau$ with $\modd\tau\leq1.1\eta$, and $\modd{\alpha-1}+\modd\tau\leq 3.1\eta$. 

Therefore, unless both functions $\equiv0$, the sum in (\ref{shift}) has size $\leq 2\eta\ell(Q)/2=\eta\ell(Q)$, and we can write the new expression as
$$\alpha(x_k-(x_Q)_k+\tau\ell(Q)),$$
where $\tau$ has the right size. 

So, with the variables $x_j$ ($j\not=k$) fixed, we have two functions of $x_k$: $g^{(Q)}(U^{(Q)}_{k}(x-x_Q))$ (considered as a function of $x_k$) and $U^{(Q)}(k,k)g^{(Q)}(U^{(Q)}_{k-1}(x-x_Q))$ (ditto), where, for each, the variable $x_k$ only occurs in the $k^{th}$ \dq slot'', and that is the only place where the two functions' expressions differ. We have called the first function $\phi(x_k)$ (considered as a function of $x_k$ only). The $k^{th}$ entry of its argument is $x_k-(x_Q)_k$. If we replace $x_k$ by $\alpha(x_k-(x_Q)_k+\tau\ell(Q))+(x_Q)_k$ (where $\alpha$ and $\tau$ mean what they do in the preceding paragraph), we get
$$\alpha(x_k-(x_Q)_k+\tau\ell(Q))+(x_Q)_k-(x_Q)_k=  \alpha(x_k-(x_Q)_k+\tau\ell(Q)),$$
which is the $k^{th}$ entry of $U^{(Q)}(k,k)g^{(Q)}(U^{(Q)}_{k-1}(x-x_Q))$'s argument. Recalling that $\alpha=U^{(Q)}(k,k)$, we see that, when the variables $x_j$ ($j\not=k$) are fixed, then, as functions of $x_k$ only, if we write $\phi(x_k):=g^{(Q)}(U^{(Q)}_{k}(x-x_Q))$, then $U^{(Q)}(k,k)g^{(Q)}(U^{(Q)}_{k-1}(x-x_Q))$ will equal $\alpha\phi(\alpha(x_k-(x_Q)_k+\tau\ell(Q))+(x_Q)_k)=\alpha\phi(\alpha(x_k-x_{I_k(Q)}+\tau\ell(I_k(Q)))+x_{I_k(Q)})$. The reader should note the subtle change. We have replaced $(x_Q)_k$ (the $k^{th}$ coordinate of $Q$'s center) with $x_{I_k(Q)}$ (the center of $Q$'s $k^{th}$ cartesian factor) and $\ell(Q)$ with $\ell(I_k(Q))$, because now we are dealing with functions of the single variable $x_k$.

I claim that, as functions of $x_k$ (the other variables being fixed), $g^{(Q)}(U^{(Q)}_{k}(x-x_Q))$ and $U^{(Q)}(k,k)g^{(Q)}(U^{(Q)}_{k-1}(x-x_Q))$ are a pair of functions like the ones considered in Lemma \ref{perturbdiff}. The only possible problem is that, for the second function, we have $\modd{\alpha-1}+\modd\tau\leq3.1\eta$ and not $\leq\eta$. Recall that $\eta\leq1/(20d)$. Then $3.1\eta\leq1/2$ (as Lemma \ref{perturbdiff} requires), and so we can apply the lemma with $3.1\eta$ in place of $\eta$, at the price of multiplying our final almost-orthogonality constant by $\sqrt{3.1}$. 

For every $Q\in\Dee_d$ and $1\leq k\leq d$, define, for $x\in\Rd$,
\begin{equation}\label{psikq}\psi^{(Q),k}(x):= g^{(Q)}(U^{(Q)}_{k}(x-x_Q))-U^{(Q)}(k,k)g^{(Q)}(U^{(Q)}_{k-1}(x-x_Q)).\end{equation}
We note that each $\psi^{(Q),k}$ is zero outside $Q$. We claim that, for each $k$, the family
\begin{equation}\label{affinefamily}\bigbrace{{\psi^{(Q),k}\over\modd{Q}^{1/2}}}_{Q\in\Dee_d}\end{equation}
is almost-orthogonal in $L^2(\Rd)$ with constant $\leq C\eta^{1/2}$, with $C$ an absolute constant. That will take care of (\ref{firstfamily}), since, for every $Q$,
$$I^{(Q)}=-\modd{L^{(Q)}}\sum_{k=1}^d \modd{U^{(Q)}_{k}}\psi^{(Q),k},$$
and the determinants are bounded by an absolute constant. The upshot will be that (\ref{firstfamily}) is almost-orthogonal with constant $\leq Cd\eta^{1/2}$, with $C$ an absolute constant.

Let us now fix $k$ (and remember: we are assuming $d>1$). Slightly abusing notation, we write $x=(x_1,\ldots,x_d):=(x_k,x_k')$, where $x_k'\in\R^{d-1}$ has all the variables except $x_k$. We note that, if $h\in L^1(\Rd)$, then (with another slight abuse of notation)
$$\int h(x)\, dx=\int h(x_k,x'_k)\,dx_k\,dx'_k.$$
If $Q=\prod_1^d I_j(Q)$ is a dyadic cube then 
$$\chi_Q(x)=\chi_Q(x_1,\ldots,x_d)=\chi_{I_k(Q)}(x_k)\cdot\prod_{j\not=k}\chi_{I_j(Q)}(x_j)=:\chi_{I_k(Q)}(x_k)\cdot F_{Q,k}(x_k').$$
Since every $\psi^{(Q),k}=0$ outside $Q$, we can write
\begin{eqnarray*}{\psi^{(Q),k}(x)\over\modd{Q}^{1/2}}&=&{\psi^{(Q),k}(x_k,x_k')\over\modd{Q}^{1/2}}\cdot\chi_{I_k(Q)}(x_k)\cdot F_{Q,k}(x_k')\\
&=&{\psi^{(Q),k}(x_k,x_k')\over\modd{I_k(Q)}^{1/2}}\cdot{\modd{I_k(Q)}^{1/2}\over\modd{Q}^{1/2}}\chi_{I_k(Q)}(x_k)\cdot F_{Q,k}(x_k'),\end{eqnarray*}
where $\modd{I_k(Q)}$ means the interval's one-dimensional Lebesgue measure.

Let $\Eff\subset\Dee_d$ be finite and $\{\lambda_Q\}_{Q\in\Eff}\subset{\bf C}$. We consider the finite linear sum
$$\sum_{Q\in\Eff}\lambda_Q{\psi^{(Q),k}(x)\over\modd{Q}^{1/2}}=
\sum_{Q\in\Eff}\lambda_Q{\psi^{(Q),k}(x_k,x_k')\over\modd{I_k(Q)}^{1/2}}\cdot{\modd{I_k(Q)}^{1/2}\over\modd{Q}^{1/2}}\chi_{I_k(Q)}(x_k)\cdot F_{Q,k}(x_k').$$
By Lemma \ref{perturbdiff}, for every fixed $x_k'\in\R^{d-1}$,
\begin{eqnarray*}& \int_\R \bigmodd{\sum_{Q\in\Eff}\lambda_Q{\psi^{(Q),k}(x_k,x_k')\over\modd{I_k(Q)}^{1/2}}\cdot{\modd{I_k(Q)}^{1/2}\over\modd{Q}^{1/2}}\chi_{I_k(Q)}(x_k)\cdot F_{Q,k}(x_k')}^2\,dx_k\\
&\leq C\eta\sum_{Q\in\Eff}\modd{\lambda_Q}^2{\modd{I_k(Q)}\over\modd{Q}}F_{Q,k}(x_k'),\end{eqnarray*}
because $F_{Q,k}(x_k')=$ 0 or 1 everywhere. We are using the fact that, for every fixed $x_k'$, the non-identically-zero members of 
$$\{\psi^{(Q),k}(x_k,x_k')\cdot F_{Q,k}(x_k')\}_{Q\in\Dee_d},$$
considered as functions of $x_k\in\R$, form a family like the set $\{g^{(I)}\}_{I\in\Dee_1}$, defined by (\ref{diffofgs}) in the statement of Lemma \ref{perturbdiff}. For every $Q$,
$$\int_{\R^{d-1}}{\modd{I_k(Q)}\over\modd{Q}}F_{Q,k}(x_k')\,dx_k'=1.$$
Therefore
\begin{eqnarray*}\int_\Rd \bigmodd{\sum_{Q\in\Eff}\lambda_Q{\psi^{(Q),k}(x)\over\modd{Q}^{1/2}}}^2\,dx&=&\int_{\R^{d-1}}\left(\int_\R\bigmodd{\sum_{Q\in\Eff}\lambda_Q{\psi^{(Q),k}(x)\over\modd{Q}^{1/2}}}^2\,dx_k\right)\,dx_k'\\
&\leq& C\eta\sum_{Q\in\Eff}\modd{\lambda_Q}^2\end{eqnarray*}
(with $C$ absolute), which was to be proved: we have controlled (\ref{firstfamily}) (the \dq$I^{(Q)}$'' term).

We treat (\ref{secondfamily}) similarly, but now we move through the matrices in the \emph{other} direction. Write
\begin{equation}\label{iiqsum}II^{(Q)}(x)=\sum_{k=1}^d \left(\modd{L^{(Q)}_{k-1}}f^{(Q)}(x_Q+L^{(Q)}_{k-1}(x-x_Q))-\modd{L^{(Q)}_{k}}f^{(Q)}(x_Q+L^{(Q)}_{k}(x-x_Q))\right),\end{equation}
where
\begin{eqnarray*}L^{(Q)}_{0}&:=&L^{(Q)}\\
L^{(Q)}_{d}&:=&I_d,\end{eqnarray*}
and, for $0<k<d$,
$$L^{(Q)}_k(i,j)=\cases{\delta(i,j)&if $\max(i,j)> d-k$;\cr
L^{(Q)}(i,j)&if $\max(i,j)\leq d-k$.\cr}$$
To picture this, think of what happened with the $U^{(Q)}_k$s, but going from bottom to top instead of top to bottom. If $M$ is a $d\times d$ matrix then the set of entries $\{M(i,j):\ \max(i,j)> d-k\}$ is the union of $M$'s \emph{last} $k$ rows (counting from the top) and \emph{last} $k$ columns (counting from the left), and  $\{M(i,j):\ \max(i,j)\leq d-k\}$ is everything else: a $(d-k)\times (d-k)$ submatrix in $M$'s upper left corner. As $k$ increases from 0 to $d$, $L^{(Q)}_k$---like $U^{(Q)}_k$ earlier---changes from $L^{(Q)}$ into $I_d$, one elbow-shaped layer at a time, with the elbow joints opening and tending toward the upper left (northwestward) instead of the lower right. 

Just as the $U^{(Q)}_k$s were all upper triangular, the $L^{(Q)}_k$s are all lower triangular. Also, for every $0\leq k\leq d$, $\nrm{I_d-L^{(Q)}_k}_\infty\leq\nrm{I_d-L^{(Q)}}_\infty$---and for the same reason.

Arguing as before, suppose that $z=\{z_j\}_1^d\in\Rd$ and $\xi=\{\xi_i\}_1^d:=L^{(Q)}_k z$ (both written as column vectors). Obviously, if $k=0$ then $\xi_i=\sum_{j=1}^d L^{(Q)}(i,j)z_j$, and if $k=d$ then $\xi_i=z_i$ $(1\leq i\leq d$). If $0<k<d$ then
$$\xi_i=\cases{z_i&if $i>d-k$;\cr \sum_{j=1}^i L^{(Q)}(i,j)z_j=\sum_{j=1}^d L^{(Q)}(i,j)z_j&if $i\leq d-k$,\cr}$$
where the equality
$$\sum_{j=1}^i L^{(Q)}(i,j)z_j=\sum_{j=1}^d L^{(Q)}(i,j)z_j$$
holds for $i\leq d-k$ because $L^{(Q)}$ is \emph{lower} triangular.

Look at one of the terms in (\ref{iiqsum}):
\begin{equation}\label{iiqterm}\modd{L^{(Q)}_{k-1}}f^{(Q)}(x_Q+L^{(Q)}_{k-1}(x-x_Q))-\modd{L^{(Q)}_{k}}f^{(Q)}(x_Q+L^{(Q)}_{k}(x-x_Q)).\end{equation}
Reasoning much as we did with the $U^{(Q)}_k$s,
$$\modd{L^{(Q)}_{k}}=\cases{\prod_{1}^{d-k} L^{(Q)}(i,i)&if $0\leq k<d$;\cr
1&if $k=d$.\cr}$$
Therefore we can rewrite (\ref{iiqterm}) as
\begin{equation}\label{iiqtermrewr}\modd{L^{(Q)}_{k}}\left(L^{(Q)}(d-k+1,d-k+1)f^{(Q)}(x_Q+L^{(Q)}_{k-1}(x-x_Q))
-f^{(Q)}(x_Q+L^{(Q)}_{k}(x-x_Q))\right).\end{equation}
The determinant at the far left of (\ref{iiqtermrewr}) is bounded by an absolute constant (same reason as with the $U^{(Q)}_k$s). As we did earlier, we will ignore it and focus on
\begin{eqnarray}\label{reducediiqtermrewr}&L^{(Q)}(d-k+1,d-k+1)f^{(Q)}(x_Q+L^{(Q)}_{k-1}(x-x_Q))\\
&-f^{(Q)}(x_Q+L^{(Q)}_{k}(x-x_Q))\nonumber.\end{eqnarray}
Once again, for ease of reading, we write $x-x_Q=z=(z_1,\ldots,z_d)$. We consider the first piece of (\ref{reducediiqtermrewr}). It is
$$L^{(Q)}(d-k+1,d-k+1))f^{(Q)}(x_Q+L^{(Q)}_{k-1}(x-x_Q)).$$
Written as a column, the entries in $f^{(Q)}$'s argument are:
\begin{eqnarray*}&(x_Q)_1+L^{(Q)}(1,1)z_1\\
&(x_Q)_2+L^{(Q)}(2,1)z_1+L^{(Q)}(2,2)z_2\\
&\vdots\\
&(x_Q)_{d-k+1}+\sum_{j=1}^{d-k+1}L^{(Q)}(d-k+1,j)z_j\\
&(x_Q)_{d-k+2}+z_{d-k+2}\qquad ({\rm i.e.,\ }x_{d-k+2})\\
&\vdots\\
&(x_Q)_d+z_d\qquad ({\rm i.e.,\ }x_{d})\end{eqnarray*}
unless $k=1$, in which case there are no entries beyond the $d-k+1=d$ slot (which, recall, means $L^{(Q)}_{k-1}=L^{(Q)}_0=L^{(Q)}$). 
The $(d-k+1)^{th}$ entry in the last expression is
$$(x_Q)_{d-k+1}+L^{(Q)}(d-k+1,d-k+1)z_{d-k+1}+\sum_{j=1}^{d-k} L^{(Q)}(d-k+1,j)z_j,$$
and it is the only entry containing the variable $z_{d-k+1}$ (or, equivalently, $x_{d-k+1}$). Let's look at the other piece of (\ref{reducediiqtermrewr}). It's
$$f^{(Q)}(x_Q+L^{(Q)}_{k}(x-x_Q)).$$
When we write the entries of its argument as a column, we get:
\begin{eqnarray*}&(x_Q)_1+L^{(Q)}(1,1)z_1\\
&(x_Q)_2+L^{(Q)}(2,1)z_1+L^{(Q)}(2,2)z_2\\
&\vdots\\
&(x_Q)_{d-k}+\sum_{j=1}^{d-k}L^{(Q)}(d-k,j)z_j\\
&(x_Q)_{d-k+1}+z_{d-k+1}\qquad ({\rm i.e.,\ }x_{d-k+1})\\
&\vdots\\
&(x_Q)_d+z_d \qquad({\rm i.e.,\ }x_{d}).\end{eqnarray*}
(Notice that, if $k=d$---meaning $L^{(Q)}_k=L^{(Q)}_d=I_d$---we get $x_1$, \dots, $x_d$ down the line, as we should!) The point is that $f^{(Q)}(x_Q+L^{(Q)}_{k-1}(x-x_Q))$ and 
$f^{(Q)}(x_Q+L^{(Q)}_{k}(x-x_Q))$ are identical except for their $(d-k+1)^{th}$ entries and, for both of them, the variable $z_{d-k+1}$ occurs only there. This implies that $f^{(Q)}(x_Q+L^{(Q)}_{k-1}(x-x_Q))$ and 
$f^{(Q)}(x_Q+L^{(Q)}_{k}(x-x_Q))$, when written as functions of $x=(x_1,\ldots,x_d)$, are also identical
except for their $(d-k+1)^{th}$ entries, and both have the variable $x_{d-k+1}$ occurring only there. Suppose we temporarily fix all the variables except $x_{d-k+1}$. Then, as a function of $x_{d-k+1}$, we can write $f^{(Q)}(x_Q+L^{(Q)}_{k}(x-x_Q))$ (which is, note, the \emph{second} term in (\ref{reducediiqtermrewr})) as $\phi(x_{d-k+1})$ for some function $\phi:\R\to{\bf C}$. Like the one we saw earlier (when considering the \dq$I^{(Q)}$'' terms), this $\phi$ has two interesting properties. It is of bounded variation, with total variation $\leq 1$; and it is zero outside of $I_{d-k+1}(Q)$, the $(d-k+1)^{th}$ cartesian factor of $Q$. The first property is an immediate consequence of $SNBV$. The second follows from $\nrm{I_d-L^{(Q)}_k}_\infty\leq\nrm{I_d-L^{(Q)}}_\infty\leq2\eta$ and the fact that $\modd{x_{d-k+1}-(x_Q)_{d-k+1}}\geq\ell(Q)/2$ implies $\nrm{x-x_Q}_\infty\geq\ell(Q)/2$ implies $\nrm{L^{(Q)}_{k}(x-x_Q)}_\infty>\ell(Q)/6$, and {therefore}
$$f^{(Q)}(x_Q+L^{(Q)}_{k}(x-x_Q))=0,$$
because $f^{(Q)}$ is zero outside the central third of $Q$. The other term in (\ref{reducediiqtermrewr}) is $L^{(Q)}(d-k+1,d-k+1)f^{(Q)}(x_Q+L^{(Q)}_{k-1}(x-x_Q))$. As a function of $x_{d-k+1}$ only (with the other variables fixed), it can be expressed in terms of the same $\phi$. It is $\alpha\phi(\alpha(x_{d-k+1}-(x_Q)_{d-k+1}+\tau\ell(Q))+(x_Q)_{d-k+1})$, where $\alpha=L^{(Q)}(d-k+1,d-k+1)$ and $\modd{\alpha-1}+\modd\tau\leq 7\eta/2$. This is very much as we saw before, but we will repeat part of it for the reader's benefit. We begin by comparing the $(d-k+1)^{th}$ entries of $f^{(Q)}(x_Q+L^{(Q)}_{k}(x-x_Q))$ and $f^{(Q)}(x_Q+L^{(Q)}_{k-1}(x-x_Q))$. The $(d-k+1)^{th}$ entry of $f^{(Q)}(x_Q+L^{(Q)}_{k}(x-x_Q))$ is $(x_Q)_{d-k+1}+z_{d-k+1}$ (i.e., $x_{d-k+1}$). The $(d-k+1)^{th}$ entry of $f^{(Q)}(x_Q+L^{(Q)}_{k-1}(x-x_Q))$ is
\begin{eqnarray*}&(x_Q)_{d-k+1}+\sum_{j=1}^{d-k+1}L^{(Q)}(d-k+1,j)z_j\\
&=(x_Q)_{d-k+1}+L^{(Q)}(d-k+1,d-k+1)(x_{d-k+1}-(x_Q)_{d-k+1})\\
&+\sum_{j=1}^{d-k}L^{(Q)}(d-k+1,j)(x_j-(x_Q)_j)\\
&=\alpha(x_{d-k+1}-(x_Q)_{d-k+1}+\tau\ell(Q))+(x_Q)_{d-k+1},\end{eqnarray*}
with
$$\tau\ell(Q)=\alpha^{-1}\sum_{j=1}^{d-k}L^{(Q)}(d-k+1,j)(x_j-(x_Q)_j).$$
Exactly as with the $U^{(Q)}_k$s, this gives the right bound for $\modd{\alpha-1}+\modd\tau$. If $\modd{x_j-(x_Q)_j}>\ell(Q)/2$ for some $1\leq j\leq d-k$ then, as we saw above, $x_Q+L^{(Q)}_{k-1}(x-x_Q)$ and $x_Q+L^{(Q)}_{k}(x-x_Q)$ both fall outside the middle third of $Q$, and $f^{(Q)}(x_Q+L^{(Q)}_{k-1}(x-x_Q))$ and $f^{(Q)}(x_Q+L^{(Q)}_{k}(x-x_Q))$ will equal 0 no matter what the value of $x_{d-k+1}$, with the obvious consequence that
$$\phi(x_{d-k+1})-\alpha\phi(\alpha(x_{d-k+1}-(x_Q)_{d-k+1})+\tau\ell(Q)+(x_Q)_{d-k+1})$$
will be identically 0 as a function of $x_{d-k+1}$. I.e., in every non-trivial case we must have $\modd{x_j-(x_Q)_j}\leq\ell(Q)/2$ for all $1\leq j\leq d-k$. The rest of the estimate, concluding with $\modd{\alpha-1}+\modd\tau\leq 3.1\eta$, is as we saw earlier. Recall that the \dq$3.1\eta$'' is good enough here because $3.1\eta\leq1/2$.

Now we proceed much as we did with the $I^{(Q)}$ terms. For every $Q\in\Dee_d$, $1\leq k\leq d$, and $x\in\Rd$, set
$$\psi^{(Q),k}_{new}(x):=L^{(Q)}(d-k+1,d-k+1)f^{(Q)}(x_Q+L^{(Q)}_{k-1}(x-x_Q))-f^{(Q)}(x_Q+L^{(Q)}_{k}(x-x_Q)).$$
We claim that, for each $k$, the new family
$$\bigbrace{{\psi^{(Q),k}_{new}\over\modd{Q}^{1/2}}}_{Q\in\Dee_d}$$
is almost-orthogonal in $L^2(\Rd)$ with constant $\leq C\eta^{1/2}$, with $C$ an absolute constant. As with (\ref{affinefamily}) and (\ref{firstfamily}), this will take care of (\ref{secondfamily}), yielding, as before, an almost-orthogonality constant $\leq C d\eta^{1/2}$ for (\ref{secondfamily}), with $C$ an absolute constant. The proof of this is essentially identical to that for (\ref{firstfamily}), the only difference being that $x_{d-k+1}$ plays the role $x_k$ did earlier. 

Write $x=(x_1,\ldots,x_d):=(x_{d-k+1},x_k')$, where $x_k'\in\R^{d-1}$ has all the variables except $x_{d-k+1}$. Clearly, if $h\in L^1(\Rd)$, then 
$$\int h\,dx=\int h(x_{d-k+1},x_k')\,dx_{d-k+1}\,dx_k',$$
more or less as we had earlier, with a similar abuse of notation. 
If $Q=\prod_1^d I_j(Q)$ is a dyadic cube then 
$$\chi_Q(x)=\chi_{I_{d-k+1}(Q)}(x_{d-k+1})\cdot\prod_{j\not=d-k+1}\chi_{I_j(Q)}(x_j)=:\chi_{I_{d-k+1}(Q)}(x_{d-k+1})\cdot G_{Q,k}(x_k').$$
Since every $\psi^{(Q),k}_{new}=0$ outside $Q$, we can write
\begin{eqnarray*}{\psi^{(Q),k}_{new}(x)\over\modd{Q}^{1/2}}&=&{\psi^{(Q),k}_{new}(x_{d-k+1},x_k')\over\modd{Q}^{1/2}}\cdot\chi_{I_{d-k+1}(Q)}(x_{d-k+1})\cdot G_{Q,k}(x_k')\\
&=&{\psi^{(Q),k_{new}}(x_{d-k+1},x_k')\over\modd{I_{d-k+1}(Q)}^{1/2}}\cdot{\modd{I_{d-k+1}(Q)}^{1/2}\over\modd{Q}^{1/2}}\chi_{I_{d-k+1}(Q)}(x_{d-k+1})\cdot G_{Q,k}(x_k'),\end{eqnarray*}
where $\modd{I_{d-k+1}(Q)}$ means the interval's one-dimensional Lebesgue measure.

Let $\Eff\subset\Dee_d$ be finite and $\{\lambda_Q\}_{Q\in\Eff}\subset{\bf C}$, and consider the finite linear sum
\begin{eqnarray*}&\sum_{Q\in\Eff}\lambda_Q{\psi^{(Q),k}_{new}(x)\over\modd{Q}^{1/2}}\\
&=\sum_{Q\in\Eff}\lambda_Q{\psi^{(Q),k}_{new}(x_{d-k+1},x_k')\over\modd{I_{d-k+1}(Q)}^{1/2}}\cdot{\modd{I_{d-k+1}(Q)}^{1/2}\over\modd{Q}^{1/2}}\chi_{I_{d-k+1}(Q)}(x_k)\cdot G_{Q,k}(x_k').\end{eqnarray*}
By Lemma \ref{perturbdiff}, for every fixed $x_k'\in\R^{d-1}$ (using, recall, the more recent definition of $x_k'$),
\begin{eqnarray*}& \int_\R \bigmodd{\sum_{Q\in\Eff}\lambda_Q{\psi^{(Q),k}_{new}(x_{d-k+1},x_k')\over\modd{I_{d-k+1}(Q)}^{1/2}}\cdot{\modd{I_{d-k+1}(Q)}^{1/2}\over\modd{Q}^{1/2}}\chi_{I_{d-k+1}(Q)}(x_{d-k+1})\cdot G_{Q,k}(x_k')}^2\,dx_{d-k+1}\\
&\leq C\eta\sum_{Q\in\Eff}\modd{\lambda_Q}^2{\modd{I_{d-k+1}(Q)}\over\modd{Q}}G_{Q,k}(x_k'),\end{eqnarray*}
because $G_{Q,k}(x_k')=$ 0 or 1 everywhere. For every $Q$,
$$\int_{\R^{d-1}}{\modd{I_{d-k+1}(Q)}\over\modd{Q}}G_{Q,k}(x_k')\,dx_k'=1.$$
Therefore
\begin{eqnarray*}\int_\Rd \bigmodd{\sum_{Q\in\Eff}\lambda_Q{\psi^{(Q),k}_{new}(x)\over\modd{Q}^{1/2}}}^2\,dx&=&\int_{\R^{d-1}}\left(\int_\R\bigmodd{\sum_{Q\in\Eff}\lambda_Q{\psi^{(Q),k}_{new}(x)\over\modd{Q}^{1/2}}}^2\,dx_{d-k+1}\right)\,dx_k'\\
&\leq& C\eta\sum_{Q\in\Eff}\modd{\lambda_Q}^2,\end{eqnarray*}
which was to be proved: the \dq$II^{(Q)}$'' family also has the right bound.

Now we remove our two preliminary assumptions, in the reverse of the order in which we imposed them. The second assumption was that every $y^{(Q)}=0$, so we treat it first. 

We have been using $\widetilde{f^{(Q)}}$ to mean 
$$ f^{(Q)}(x_Q+A^{(Q)}(x-x_Q)).$$
We now reassign it to its original meaning from (\ref{originalperturb}):
$$ \widetilde{f^{(Q)}}(x):=f^{(Q)}(x_Q+A^{(Q)}(\ell(Q)y^{(Q)}+x-x_Q)).$$
We will use
$$\widetilde{f_0^{(Q)}}(x):=f^{(Q)}(x_Q+A^{(Q)}(x-x_Q))$$
to stand for the \dq$y^{(Q)}=0$'' special case we have been considering. 

We claim that
$$\bigbrace{{\widetilde{f^{(Q)}}-\widetilde{f_0^{(Q)}}\over\modd{Q}^{1/2}}}_{Q\in\Dee_d}$$
is almost-orthogonal with constant $\leq Cd\eta^{1/2}$, with $C$ an absolute constant. We prove the claim \dq one step at a time'', more or less the way we handled the $I^{(Q)}$ and $II^{(Q)}$ families above. Define $y^{(Q)}=:(y^{(Q)}_1,\,y^{(Q)}_2,\ldots,y^{(Q)}_d)$. For $0\leq k\leq d$ define
\begin{eqnarray*}w^{(Q)}_0&:=&0\\
w^{(Q)}_1&:=&(y^{(Q)}_1,0,0,\ldots,0)\\
w^{(Q)}_2&:=&(y^{(Q)}_1,y^{(Q)}_2,0,0,\ldots,0)\\
&\vdots\\
w^{(Q)}_d&:=&(y^{(Q)}_1,\ldots,y^{(Q)}_d)=y^{(Q)},\end{eqnarray*}
and similarly set
$$\widetilde{f_k^{(Q)}}(x):=f^{(Q)}(x_Q+A^{(Q)}(x-x_Q+\ell(Q)w^{(Q)}_k)).$$

In the following it will be helpful to remember that $w_k^{(Q)}$ is always a vector, $y^{(Q)}_k$ is always a number, and $w^{(Q)}_k$'s components are taken from 
$$\{0,y^{(Q)}_1,y^{(Q)}_2,y^{(Q)}_3,\ldots,y^{(Q)}_d\}.$$
(The author had trouble keeping them straight.)

We note that $\widetilde{f_d^{(Q)}}=\widetilde{f^{(Q)}}$ and that, as $k$ goes from 0 to $d$, $\widetilde{f_k^{(Q)}}$ changes from $\widetilde{f_0^{(Q)}}$ to $\widetilde{f^{(Q)}}$, one variable at a time. If $1\leq k\leq d$ then, as functions of $x=(x_1,\ldots,x_d)$, $\widetilde{f_k^{(Q)}}$ and $\widetilde{f_{k-1}^{(Q)}}$ differ only in their $k^{th}$ entries, and then by a constant with absolute value $\leq \eta\ell(Q)$. 

We write
$$\widetilde{f^{(Q)}}-\widetilde{f_0^{(Q)}}=\widetilde{f_d^{(Q)}}-\widetilde{f_0^{(Q)}}=\sum_{k=1}^d \left(\widetilde{f_k^{(Q)}}-\widetilde{f_{k-1}^{(Q)}}\right).$$
We shall show that, for each $1\leq k\leq d$, the family
$$\bigbrace{{\widetilde{f_k^{(Q)}}-\widetilde{f_{k-1}^{(Q)}}\over\modd{Q}^{1/2}}}_{Q\in\Dee_d}$$
is almost-orthogonal with constant $\leq C\eta^{1/2}$ ($C$ absolute). This part will go fast.

As we saw before (and by the same argument), every $\widetilde{f_k^{(Q)}}$ will belong to $SNBV(Q)$.  Specifically, every $\widetilde{f_k^{(Q)}}$ is 0 outside $Q$. Let's temporarily fix $1\leq k\leq d$ and define $\phi^{(Q),k}(x):=\widetilde{f_{k-1}^{(Q)}}(x)$. As we did before, write $x\in\Rd$ as $(x_k,x_k')$, where $x_k'\in\R^{d-1}$ consists of all the other variables. Then $\widetilde{f_{k}^{(Q)}}(x)=\widetilde{f_{k}^{(Q)}}(x_k,x_k')$ is simply $\phi^{(Q),k}(x_k+\ell(Q)y^{(Q)}_k,x_k')$. We have reduced our problem to showing that, for each $1\leq k\leq d$, the family
$$\bigbrace{{\phi^{(Q),k}(x_k+\ell(Q)y^{(Q)}_k,x_k')-\phi^{(Q),k}(x_k,x_k')}\over\modd{Q}^{1/2}}_{Q\in\Dee_d}$$
is almost-orthogonal with constant $\leq C\eta^{1/2}$. But we've already done this. Much as we did in (\ref{psikq}), we define 
$$\psi^{(Q),k}(x_k,x_k'):=\phi^{(Q),k}(x_k+\ell(Q)y^{(Q)}_k,x_k')-\phi^{(Q),k}(x_k,x_k').$$
We seek to show that
$$\bigbrace{{\psi^{(Q),k}\over\modd{Q}^{1/2}}}_{Q\in\Dee_d}$$
is almost-orthogonal in $L^2(\Rd)$ with constant $\leq C\eta^{1/2}$. This is just the situation with (\ref{affinefamily}) all over again, and even easier, because now the only difference between $\phi^{(Q),k}(x_k+\ell(Q)y^{(Q)}_k,x_k')$ and $\phi^{(Q),k}(x_k,x_k')$ is a small translation in the $x_k$-variable. We refer the reader back to the two full paragraphs following (\ref{affinefamily}), starting with \dq Let us now fix $k$ (and remember: we are assuming $d>1$)'' and ending with \dq \dots\ we have controlled (\ref{firstfamily}) (the \dq$I^{(Q)}$'' term)'', for the detailed argument.

Finally we remove the assumption that every $f^{(Q)}=0$ outside the central third of $Q$. We use an elementary result from \cite{Wilsonbook} (Theorem 5.3, p.~91). For every $Q\in\Dee_d$ let $\widetilde Q$ mean its concentric triple and set $\widetilde{\Dee_d}:=\{\widetilde Q:\ Q\in\Dee_d\}$. Theorem 5.3 from \cite{Wilsonbook} states that we can write $\widetilde{\Dee_d}$ as a union of $3^d$ pairwise disjoint families $\{\Gee_k\}_1^{3^d}$, where each $\Gee_k$ has certain good properties. These properties are:\smallskip

i) $\forall Q,Q'\in\Gee_k$ we have $Q\subset Q'$, $Q'\subset Q$, or $Q\cap Q'=\emptyset$.\par
ii) Every $Q\in\Gee_k$ is the union of $2^d$ cubes from $\Gee_k$, each having sidelength equal to $(1/2)\ell(Q)$.\par
iii) Every $Q\in\Gee_k$ is contained in some $Q'\in\Gee_k$ having sidelength equal to $2\ell(Q)$.\par
iv) For every integer $j$, 
$$\Rd=\bigcup_{Q\in\Gee_k\atop \ell(Q)=3\cdot 2^j}Q.$$\smallskip 

Now we follow the procedure from \cite{WilsonBvconvao}. For each $1\leq k\leq 3^d$ set $G_k:=\{Q\in\Dee_d:\ \widetilde Q\in\Gee_k\}$. We show that, for each $k$, the family
\begin{equation}\label{triplefamily}\bigbrace{{f^{(Q)}-\modd{A^{(Q)}}\widetilde{f^{(Q)}}\over\modd{\widetilde Q}^{1/2}}}_{Q\in G_k}\end{equation}
(don't miss the tilde on the $Q$ in the denominator) is almost-orthogonal with constant $\leq Cd\eta^{1/2}$. If $Q\in G_k$ then $f^{(Q)}=0$ outside the central third of $\widetilde Q$, and the cubes of $\{\widetilde Q\}_{Q\in G_k}=\Gee_k$ mimic the essential properties of the dyadic cubes $\Dee_d$. As was done in \cite{WilsonBvconvao}, we repeat the proof of Theorem \ref{covid}, adapted to $\Gee_k$ and to triples of dyadic intervals in \R\ (the Schur test argument). We refer the reader to \cite{WilsonBvconvao} for the details. We have proved Theorem \ref{covid}. \cs\medskip

{\bf Proof of Corollary \ref{covidcor}.} With our choice of $\widetilde c(d)$, we have $\modd{1-\modd{A^{(Q)}}}<1/2$ for all $Q$. Therefore it suffices to show that
$$\bignrm{\bigbrace{{\modd{A^{(Q)}}(f^{(Q)}-\widetilde{f^{(Q)}})\over\modd{Q}^{1/2}}}_{Q\in\Dee_d}}_{AO(\Dee_d)}\leq C(d)\eta^{1/2}.$$
We also have $\modd{1-\modd{A^{(Q)}}}\leq C(d)\eta$. Now we write
$$\modd{A^{(Q)}}(f^{(Q)}-\widetilde{f^{(Q)}})=(\modd{A^{(Q)}}-1)f^{(Q)}+(f^{(Q)}-\modd{A^{(Q)}}\widetilde{f^{(Q)}}).$$
Since every $f^{(Q)}\in SNBV_0(Q)\subset NBV_0(Q)$, Theorem \ref{bvconvao} implies that
$$\bignrm{\bigbrace{{(\modd{A^{(Q)}}-1)f^{(Q)}\over\modd{Q}^{1/2}}}_{Q\in\Dee_d}}_{AO(\Dee_d)}\leq C(d)\eta\leq C(d)\eta^{1/2};$$
while
$$\bignrm{\bigbrace{{(f^{(Q)}-\modd{A^{(Q)}}\widetilde{f^{(Q)}}\over\modd{Q}^{1/2}}}_{Q\in\Dee_d}}_{AO(\Dee_d)}\leq C(d)\eta^{1/2}$$
follows from Theorem \ref{covid}. That proves the corollary. \cs \medskip 

Our final result is the proof of Corollary \ref{covidcordyadav}. To give it some context, we first prove an elementary result about averages. For $Q\in\Dee_d$, let
$Q^*$ be as defined in Corollary \ref{covidcordyadav}'s hypotheses,
$$\chi_{Q^*}(x):=\chi_Q(x_Q+A^{(Q)}(\ell(Q)y^{(Q)}+x-x_Q)),$$
where $y^{(Q)}$ and $A^{(Q)}$ are as in Theorem \ref{covid}. We think of $Q^*$ as a slight distortion of $Q$. If the $\eta$ from Theorem \ref{covid} is chosen small enough then geometric arguments like those already seen show that $(1-c_1\eta)Q\subset Q^*\subset (1+c_1\eta)Q$ for some small positive $c_1$, implying that $\modd{Q\Delta Q^*}\leq c\eta\modd Q$. Therefore, for any $g\in L^2(\Rd)$,
\begin{eqnarray*}\bigmodd{g_Q-g_{Q^*}}&\leq&\bigmodd{{1\over\modd Q}\left(\int_Q g\,dt-\int_{Q^*}g\,dt\right)}+\bigmodd{{1\over\modd Q}-{1\over\modd{Q^*}}}\int_{Q^*}\modd g\,dt\\
&\leq& {1\over\modd Q}\int_{Q\Delta Q^*}\modd g\,dt+{C\eta\over\modd{Q^*}}\int_{Q^*}\modd g\,dt,\end{eqnarray*}
where $g_Q$ and $g_{Q^*}$ mean the usual averages over $Q$ and $Q^*$. If $x\in Q$ then the second term is clearly bounded by $C\eta Mg(x)$ (the Hardy-Littlewood maximal function of $g$), while, by H\"older's inequality, for any $r>0$,
$${1\over\modd Q}\int_{Q\Delta Q^*}\modd g\,dt\leq \left({1\over\modd Q}\int_{Q\cup Q^*} \modd g^r\,dt\right)^{1/r}\left({\modd{Q\Delta Q^*}\over\modd Q}\right)^{1/r'}\leq C\eta^{1/r'}M_rg(x),$$
where $M_r(\cdot)$ is the $r^{th}$ power Hardy-Littlewood maximal function,
$$M_rg(x):=\sup_{x\in Q\atop Q{\rm\ a\ cube}}\left({1\over\modd Q}\int_Q\modd g^r\,dt\right)^{1/r}$$
and $r'=r$'s dual index. Thus, for $0<r<2$,
$$\sup_{x\in Q\in\Dee_d}\bigmodd{g_Q-g_{Q^*}}\leq C(d,r)\eta^{1/r'}M_r(g)(x).$$
$M_r(\cdot)$ is bounded on $L^2$ for $r<2$, so we obtain
\begin{equation}\label{maxfunceta}\bignrm{\sup_{x\in Q\in\Dee_d}\bigmodd{g_Q-g_{Q^*}}}_2\leq C(d,r)\eta^{1/r'}\nrm g_2,\end{equation}
where the exponent on $\eta$---$1/r'$---is strictly smaller than $1/2$.
\medskip

{\bf Proof of Corollary \ref{covidcordyadav}.} Apply Theorem \ref{covid} to $f^{(Q)}:=\chi_Q$. Then $\widetilde{f^{(Q)}}=\chi_{Q^*}$, and $\modd{A^{(Q)}}\widetilde{f^{(Q)}}={\modd Q\over\modd{Q^*}}\chi_{Q^*}$. For any $g\in L^2(\Rd)$,
\begin{equation}\label{maxsumineq}\sum_{Q\in\Dee_d}\bigmodd{\left\langle g,{f^{(Q)}-\modd{A^{(Q)}}\widetilde{f^{(Q)}}\over\modd Q^{1/2}}\right\rangle}^2\leq C(d)\eta\nrm g_2^2.\end{equation}
For any $Q\in\Dee_d$,
$$\left\langle g,{f^{(Q)}-\modd{A^{(Q)}}\widetilde{f^{(Q)}}\over\modd Q^{1/2}}\right\rangle=\modd Q^{1/2}\left(g_Q-g_{Q^*}\right).$$
We can rewrite the left-hand side of (\ref{maxsumineq}) as
$$\sum_{Q\in\Dee_d}\modd Q\bigmodd{g_Q-g_{Q^*}}^2=\int_\Rd\left(\sum_{Q\in\Dee_d\atop x\in Q}\bigmodd{g_Q-g_{Q^*}}^2\right)\,dx$$
and obtain
$$\bignrm{\left(\sum_{Q\in\Dee_d}\bigmodd{g_Q-g_{Q^*}}^2\chi_Q\right)^{1/2}}_2\leq C(d)\eta^{1/2}\nrm g_2.$$
To see that the exponent is sharp, put $d=1$, $g:=\chi_{[0,\eta)}$ for $\eta$ small, and $Q^*:=Q$ for every $Q\in\Dee_1$ except for $[0,1)$, with $[0,1)^*:=[0,1)+\eta$. Then $\nrm{g}_2=\eta^{1/2}$, $g_{[0,1)}-g_{[0,1)^*}=\eta$, and 
$$\bignrm{\left(\sum_{Q\in\Dee_1}\bigmodd{g_Q-g_{Q^*}}^2\chi_Q\right)^{1/2}}_2=\left(\int \modd{g_{[0,1)}-g_{[0,1)^*}}^2\chi_{[0,1)}\,dx\right)^{1/2}=\eta=\eta^{1/2}\nrm{g}_2.$$
\cs

\noindent Department of Mathematics\\
University of Vermont\\
Burlington, Vermont\ \ 05405\\
USA\\
Email address: Mike.Wilson@uvm.edu

\end{document}